\documentclass[11pt]{amsart}
\usepackage[left=3cm,top=3cm,right=3cm,bottom=3cm]{geometry}
\usepackage{amsfonts}
\usepackage{enumerate,subfigure,color}
\usepackage{amssymb,latexsym,amsmath,graphicx}
\usepackage{framed}
\usepackage{afterpage}

\newcommand{\dbar}{\bar{\partial}}

\newcommand{\Om}{\Omega}

\newcommand{\R}{{\mathbb R}}

\newcommand{\C}{{\mathbb C}}

\DeclareMathOperator{\de}{\partial}

\DeclareMathOperator{\dez}{\de_z}
\DeclareMathOperator{\dbarz}{\dbar_z}

\DeclareMathOperator{\by}{\times}
\DeclareMathOperator{\bndry}{\partial\Omega}
\DeclareMathOperator{\gprior}{\gamma_{\mbox{\tiny{\textbf{PR}}}}}
\DeclareMathOperator{\prior}{\mbox{\tiny{\textbf{PR}}}}
\DeclareMathOperator{\New}{\mbox{\tiny{\textbf{new}}}}
\DeclareMathOperator{\Mint}{M_{{ij}}^{\mbox{\tiny{int}}}}
\DeclareMathOperator{\Minta}{M_{{11}}^{\mbox{\tiny{int}}}}
\DeclareMathOperator{\Mintb}{M_{{12}}^{\mbox{\tiny{int}}}}
\DeclareMathOperator{\Mintc}{M_{{21}}^{\mbox{\tiny{int}}}}
\DeclareMathOperator{\Mintd}{M_{{22}}^{\mbox{\tiny{int}}}}

\DeclareMathOperator{\gDB}{\gamma_{\mbox{\tiny{\textbf{DB}}}}}
\DeclareMathOperator{\Mprior}{M^{\mbox{\tiny{\textbf{PR}}}}}
\DeclareMathOperator{\Qprior}{Q^{\mbox{\tiny{\textbf{PR}}}}}
\DeclareMathOperator{\Sprior}{S^{\mbox{\tiny{\textbf{PR}}}}}

\DeclareMathOperator{\gnew}{\gamma_{\New}}

\date{February 23, 2016}

\begin{document}
\title[Admittivity D-bar EIT with a Spatial Prior]{Incorporating a Spatial Prior into Nonlinear D-Bar EIT imaging for Complex Admittivities}
\author{S.~J.~Hamilton}\thanks{S.~J.~Hamilton is with the Department of Mathematics, Statistics, and Computer Science; Marquette University, Milwaukee, WI, 53233 USA (e-mail: sarah.hamilton@marquette.edu).}
\author{J.~L.~Mueller}\thanks{J. L. Mueller is with the Department of Mathematics and School of Biomedical Engineering, Colorado State University, CO 80523 USA (e-mail: mueller@math.colostate.edu).}
\author{M.~Alsaker}\thanks{M.~Alsaker is with Department of Mathematics; Colorado State University, Fort Collins, C0, 80523 USA (e-mail: alsaker@math.colostate.edu).}


\begin{abstract}
Electrical Impedance Tomography (EIT) aims to recover the internal conductivity and permittivity distributions of a body from electrical measurements taken on electrodes on the surface of the body.  The reconstruction task is a severely ill-posed nonlinear inverse problem that is highly sensitive to measurement noise and modeling errors.  Regularized D-bar methods have shown great promise in producing noise-robust algorithms by employing a low-pass filtering of nonlinear (nonphysical) Fourier transform data specific to the EIT problem. Including prior data with the approximate locations of major organ boundaries in the scattering transform provides a means of extending the radius of the low-pass filter to include higher frequency components in the reconstruction, in particular, features that are known with high confidence.  This information is additionally included in the system of D-bar equations with an independent regularization parameter from that of the extended scattering transform.  In this paper, this approach is used in the 2-D D-bar method for admittivity (conductivity as well as permittivity) EIT imaging.  Noise-robust reconstructions are presented for simulated EIT data on chest-shaped phantoms with a simulated pneumothorax and pleural effusion.  No assumption of the pathology is used in the construction of the prior, yet the method still produces significant enhancements of the underlying pathology (pneumothorax or pleural effusion) even in the presence of strong noise.
\end{abstract}

\maketitle

\section{Introduction}
%
%
%
%
Electrical Impedance Tomography (EIT) is a non-invasive radiation-free imaging modality in which low amplitude current is applied through electrodes placed on the surface of a body and the resulting voltages are measured.  From these surface measurements, images of the interior conductivity and permittivity can be obtained.   The severe ill-posedness of the inverse conductivity/permittivity problem limits the spatial resolution of the reconstructed images, which hinders their clinical applicability.  The use of spatial {\em a priori} information in the solution of the inverse problem provides a means of including anatomical information that is present with high confidence, while still allowing unknown features such as lung pathologies to emerge in the reconstructed image without any assumption of their presence.  In patients with serious respiratory illness, it is often the case that a CT scan is performed to obtain a diagnosis or for a regular exam in the case of a chronic illness, and the condition is monitored with one or more follow-up scans.  The initial scan can provide basic {\em a priori} information for the reconstruction algorithm such as chest shape, and approximate lung and heart sizes, and relative positions in the plane of the electrodes.  

{\em A priori} information has been used successfully in iterative reconstruction algorithms to enhance image quality  \cite{Avis95, Baysal98, CamargoThesis, Dehghani99, Dobson1994,  Ferrario12, FloresTapia10, Soleimani06, Vauhkonen1998a}, and more recently in \cite{AlsakerMueller2015} in the direct 2-D D-bar method for (real-valued) conductivity reconstruction.  In this paper, the method of \cite{AlsakerMueller2015} is extended to the 2-D D-bar algorithm for the reconstruction of complex admittivities \cite{Hamilton_Thesis_2012,Hamilton2012,Hamilton2013}.  The reconstruction algorithm for complex admittivities differs from the D-bar algorithm for real-valued conductivities in the construction of the complex geometrical optics (CGO) solutions.  While the well-developed real-valued case \cite{DeAngelo2010,Dodd2014,Isaacson2004,Isaacson2006, Mueller2003,Knudsen2007,Murphy2009,Siltanen2000} utilizes the familiar transformation of the generalized Laplace equation governing the physical EIT problem to a Schr\"odinger equation, the complex admittivity algorithm requires transforming the problem to a first order elliptic system and constructing two sets of CGO solutions.   The algorithm is described briefly in Section 2, and the reader is referred to \cite{Hamilton_Thesis_2012,Hamilton2012,Hamilton2013} for further detail.
 
The method incorporates spatial \textit{a priori} information about the admittivity distribution in the scattering transform, as well as in the system of D-bar equations, and includes regularization parameters in each place that can be adjusted to control the amount of influence the prior has on the reconstruction. The effectiveness of the method is tested here on simulated data with 0.1\% and 1.0\% added Gaussian relative noise for a 2-D phantom chest with a simulated pleural effusion and with a simulated pneumothorax.  No \textit{a priori} information about the presence of the effusion or the pneumothorax is used in the reconstruction, only \textit{a priori} spatial information about the heart and lung boundaries.  Nevertheless, both the effusion and pneumothorax become considerably sharper than in images computed without the \textit{a priori} organ boundary information.  

\begin{figure}[!h]
\centering
{\includegraphics[width=100pt]{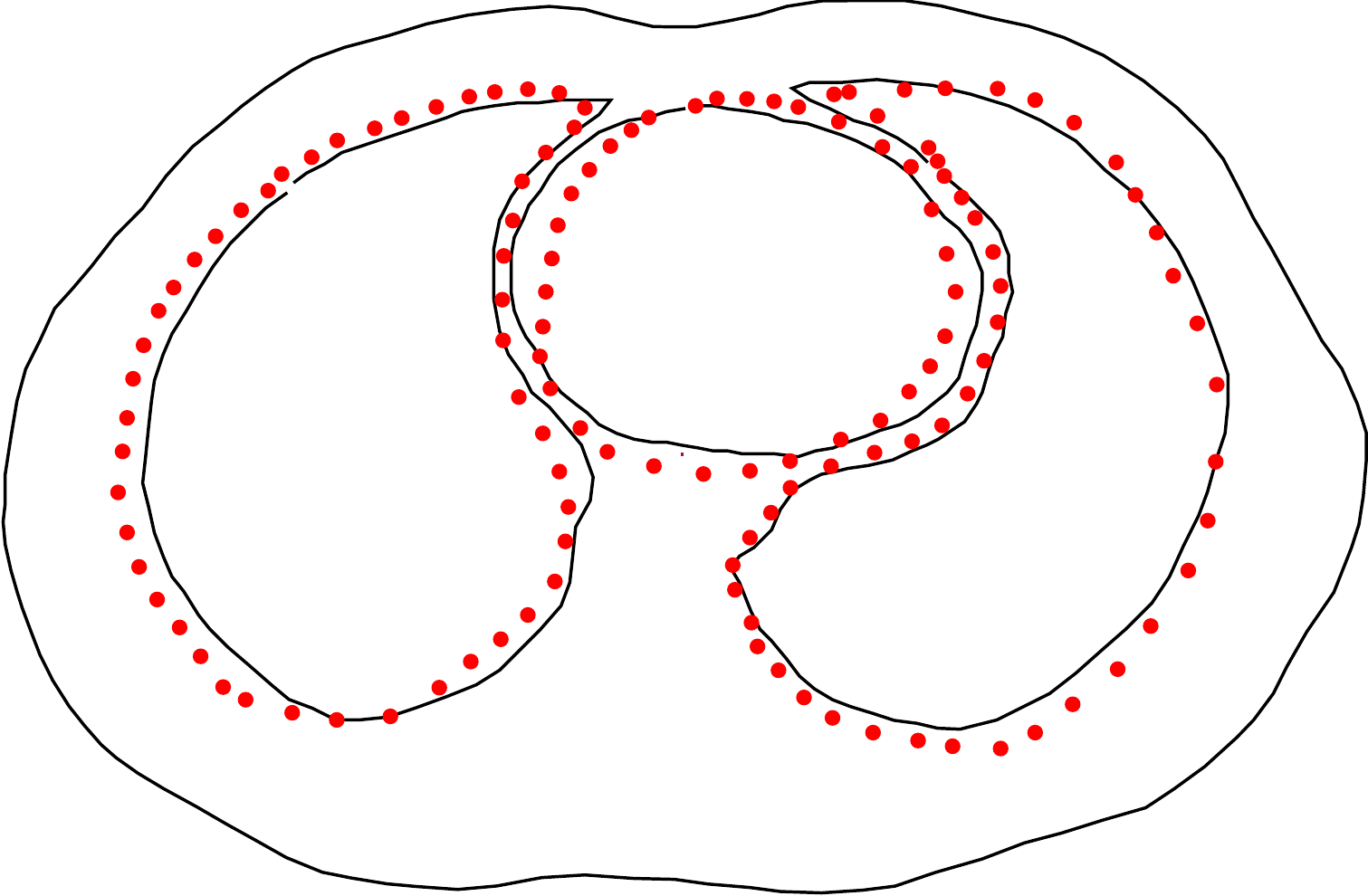}}
\caption{\label{Fig:heartnlung_prior} The ``heart and lungs'' phantom.  The true organ boundaries are depicted by the black lines, whereas the red dots depict the organ boundaries used in the prior. }
\end{figure}

The ``heart and lungs prior'' is depicted in Figure \ref{Fig:heartnlung_prior}.   While the initial prior is piecewise constant, after conductivity and permittivity values have been assigned, the prior is mollified to obtain a smooth function since the method of computing the scattering transform for the prior requires that it be differentiated.   Assigning the initial admittivity values to the prior can be done in a number of ways, and the {\em a priori} reconstruction algorithm presented here is valid for any assignment method.  In our tests, we computed an initial reconstruction with no prior from the noisy data (which we will refer to as a standard D-bar reconstruction), then computed the average conductivity/permittivity in each region of the piecewise constant ``heart and lungs prior'' prior and assigned those values to each region of the piecewise constant prior.  Further implementation details are found in Section \ref{sec:implementation}.

The paper is organized as follows.  The \textit{a priori} method is presented in Section~\ref{sec:Methods}, which first provides a brief description of the forward model in Subsection~\ref{sec:FP} used to simulate the EIT data, followed by a summary of the D-bar method for complex admittivity imaging in Subsection~\ref{sec:Dbar}, with the modifications for the \textit{a priori} method described in Subsection~\ref{sec:a_priori}.  The D-bar method for admittivity reconstruction is admittedly mathematically complicated, and the reader is referred to the papers \cite{Hamilton2013, Hamilton2012, Herrera2015} for further details.   The numerical implementation of the method is outlined in Section~\ref{sec:implementation}, the test problems described in Section~\ref{sec:Examples}, and the discussion and conclusions presented in Section~\ref{Sec:discussion}.

%


\section{Methods}\label{sec:Methods}

\subsection{The forward model} \label{sec:FP}

The electric potential $u(x,y)$ inside the 2-D region $\Omega$ is modeled by the  \emph{admittivity equation}, a generalized Laplace equation,
\begin{equation}\label{eq:admitt_eq}
\nabla\cdot\gamma(x,y)\nabla u(x,y)=0,
\end{equation}
where $\gamma(x,y)=\sigma(x,y)+i\omega\epsilon(x,y)$ denotes the complex valued admittivity, $\sigma(x,y)$ the electrical conductivity (bounded away from zero $0<\sigma(x,y)<C$), $\epsilon(x,y)$ the electrical permittivity (assumed to be non-negative), and $\omega$  the angular frequency of the applied current.  The boundary data for the inverse problem is the \textit{Dirichlet-to-Neumann} (DN) map $\Lambda_\gamma$ which maps a voltage at the boundary to the corresponding current density, i.e.,
\begin{equation}\label{eq:DNmap}
\Lambda_\gamma:u|_{\bndry}\mapsto\left.\gamma\frac{\partial u}{\partial\nu}\right|_{\bndry},
\end{equation}
where $\nu$ denotes the outward unit normal vector to the boundary $\bndry$.  In practice, to dampen rather than amplify the noise in the measured data, currents are applied and the resulting voltages are measured.  This corresponds to knowledge of the \textit{Neumann-to-Dirichlet} ND map
\begin{equation}\label{eq:NDmap}
R_\gamma:\left.\gamma\frac{\partial u}{\partial\nu}\right|_{\bndry}\mapsto u|_{\bndry}.
\end{equation}
Ensuring conservation of charge and specifying a ground, the ND map can be inverted to obtain the DN map $\Lambda_\gamma=\left(R_\gamma\right)^{-1}$.

For the simulation of the data, a finite element implementation of the complete electrode model (CEM) was used.  The CEM \cite{Cheng1989}  takes into account both the shunting effect of
the electrodes and the contact impedances between the electrodes and
tissue.  The complete electrode model consists of the admittivity
equation (\ref{eq:admitt_eq}) and the following boundary
conditions on $L$ electrodes:
$$\begin{array}{rl}
u + z_l\gamma\frac{\partial u}{\partial \nu} = U_l, & x\in e_l,
\mbox{ }l = 1,2,...,L \\\\\int_{e_l}\gamma\frac{\partial u}{\partial
\nu}dS = J_l, & x\in e_l, \mbox{ }l = 1,2,...,L
\\\\\gamma\frac{\partial u}{\partial \nu} = 0, & x\in \partial\Omega
\backslash \cup^L_{l=1}e_l,
\end{array}$$
where $z_l$ is the effective contact impedance between the $l^{th}$
electrode $e_l$ and the medium, $J_l$ is the applied current, and $U_l$ is the measured voltage.  In addition, Kirchhoff's Law and the choice of ground must be imposed to ensure
existence and uniqueness of the result:
$$\begin{array}{ccc}
\sum^L_{l=1} J_l =0, & \mbox{ and } & \sum^L_{l=1} U_l =0.
\end{array} \label{voltage_condition}$$
The uniqueness and existence of a solution to the CEM has been proven in
\cite{Somersalo1992}.

\subsection{The D-bar method for complex admittivities}\label{sec:Dbar}
D-bar methods are named for the partial derivatives with respect to the complex conjugates that arise in the equations in  the methods.  The $\dbar$ operator with respect to the complex variable $z = x+iy$  and the related operator $\partial_z$ are defined by
$$\dbar_z = \frac{1}{2}\left(\frac{\partial}{\partial x} + i \frac{\partial}{\partial y}\right), ~~~ \partial_z = \frac{1}{2}\left(\frac{\partial}{\partial x} - i \frac{\partial}{\partial y}\right).$$
  Throughout the paper, $\R^2$ is associated with $\C$ via $z=(x,y)\mapsto x+iy$.

The method described below is based on the uniqueness proof for the inverse admittivity problem \cite{Francini2000}, which was completed as a constructive proof in \cite{Hamilton2012, Hamilton_Thesis_2012}.  With the introduction of a non-physical complex parameter $k$, the admittivity equation \eqref{eq:admitt_eq} admits solutions with special exponentially growing behavior known as CGO solutions.  In particular, it was shown in  \cite{Hamilton2012} that there exist separate solutions $u_1(z,k)$ and $u_2(z,k)$ to \eqref{eq:admitt_eq} such that $u_1(z,k)\sim \frac{e^{ikz}}{ik}$ and $u_2(z,k)\sim -\frac{e^{-ik\bar{z}}}{ik}$.

Defining an operator vector $\mathcal{D} = \gamma^{1/2}(\dez ,\dbarz )^T$, the change of variables 
\begin{equation}\label{eq:Q}
Q(z)=\left[\begin{array}{cc}
0 & -\frac{1}{2} \dez \log \gamma(z)\\
-\frac{1}{2} \dbarz \log \gamma(z) & 0\\
\end{array}\right],
\end{equation}
and 
$$(M_{11},M_{21})^T=e^{-ikz}\mathcal{D}u_1, ~~~ (M_{12},M_{22})^T=e^{ik\bar{z}}\mathcal{D}u_2$$
transform the admittivity equation into the first order elliptic system \cite{Francini2000}
\begin{equation}\label{eq:DkM=QM}
D_k M(z,k)-Q(z)M(z,k)=0,
\end{equation}
where 
\[\footnotesize D_k M(z,k)=\left[\begin{array}{cc}
\dbarz & 0 \\
0 & \dez
\end{array}\right]M-ik\left[\begin{array}{cc}
1 & 0\\
0 & -1
\end{array}\right] \left[\begin{array}{cc}
0 & M_{12}\\
M_{21} & 0
\end{array}\right].\]
Equation \eqref{eq:DkM=QM} has a unique solution $M(\cdot,k)$ for $M(\cdot,k)-I\in L^p(\R^2)$ for some $p>2$.

D-bar methods follow the basic computational outline:
\begin{center}
\textbf{\small DN map} $\mapsto\begin{array}{c}
\textbf{{\small Scattering}}\\
\textbf{{\small Data}}
\end{array}\mapsto\begin{array}{c}
\textbf{{\small CGO}}\\
\textbf{{\small Solutions}}
\end{array}\mapsto$ \textbf{\small Admittivity}.
\end{center} 

The scattering data is a $2\by 2$ matrix function $S(k)$, not physically measurable from the data, with zero entries on the diagonal and off-diagonal entries given by 
\begin{equation}\label{eq:ScatDom}
\begin{array}{rcl}
S_{12}(k)&=&\frac{i}{\pi}\int_{\Om} Q_{12}(z)e(z,-\bar{k})M_{22}(z,k)\;dxdy\\
S_{21}(k)&=&-\frac{i}{\pi}\int_{\Om} Q_{21}(z)e(z,k)M_{11}(z,k)\;dxdy
\end{array}
\end{equation}
where $e(z,k)\equiv e^{i(kz+\bar{k}\bar{z})}$ and \text{supp} $Q(z)\subseteq\overline{\Omega}$.

The DN map $\Lambda_\gamma$ uniquely determines the scattering data $S(k)$, and the scattering data uniquely determines the admittivity $\gamma(z)$ \cite{Francini2000}.  
However, the relationship between the scattering data and the DN map relies on the intermediate computation of the CGO solutions $u_1$ and $u_2$ on the boundary of $\Om$ as well as functions $\Psi_{12}(z,k)\equiv e^{-ik\bar{z}}M_{12}(z,k)$ and $\Psi_{21}(z,k)\equiv e^{ikz}M_{21}(z,k)$.   This is described in Step 1 below.

\vspace{1em}
\noindent\textbf{\underline{Step 1:} From Boundary Measurements to Scattering Data:}\\
For each $|k|\leq R$, solve the following two boundary integral equations 
\begin{equation}\label{eq:BIEu1u2}
\small\begin{array}{l}
u_1(z,k)= \frac{e^{ikz}}{ik}-\underset{\bndry}{\int}G_k(z-\zeta)\left(\delta\Lambda_\gamma\right)u_1(\zeta,k)\; ds(\zeta)\\
u_2(z,k)= -\frac{e^{-ik\bar{z}}}{ik}-\underset{\bndry}{\int}G_k(\bar{\zeta}-\bar{z})\left(\delta\Lambda_\gamma\right)u_2(\zeta,k)\; ds(\zeta).
\end{array}
\end{equation}
for the traces of the CGO solutions $u_1$ and $u_2$ on the boundary.  Here $G_k(z)$ denotes the Faddeev Green's function for the Laplace operator given by (see \cite{Faddeev1966,Mueller2012}),
\[G_k(z)=e^{ikz}\int_{\R^2}\frac{e^{iz\cdot \xi}}{\xi(\bar{\xi}+2k)}\;d\xi,\]
and $\delta\Lambda_\gamma=\Lambda_\gamma - \Lambda_1$ where $\Lambda_1$ denotes the DN map corresponding to a constant admittivity $\gamma=1$.

Next, compute the traces of the CGO solutions $\Psi_{12}$ and $\Psi_{21}$ from the second set of boundary integral equations 
\begin{equation}\label{eq:PsiBIEs}
\small \begin{array}{rcl}
\Psi_{12}(z,k)&=& p.v.\underset{\bndry}{\int} \frac{e^{i\bar{k}(z-\zeta)}}{4\pi(z-\zeta)}\left(\delta\Lambda_\gamma\right)u_2(\zeta,k) ds(\zeta)\\
\Psi_{21}(z,k)&=& p.v.\underset{\bndry}{\int} \overline{\left[\frac{e^{i{k}(z-\zeta)}}{4\pi(z-\zeta)}\right]}\left(\delta\Lambda_\gamma\right)u_1(\zeta,k) ds(\zeta),
\end{array}
\end{equation}
where $p.v.$ denotes the principal value of the integral.

Then, compute the scattering transforms $S_{12}(k)$ and $S_{21}(k)$:
\begin{equation}\label{eq:scatBIEs}
\begin{array}{rcr}
S_{12}(z,k)&=&
\frac{i}{2\pi}\underset{\bndry}{\int} e^{-i\bar{k}z}\Psi_{12}(z,k)\nu(z)\;ds(z)  \\
S_{21}(z,k)&=& 
-\frac{i}{2\pi}\underset{\bndry}{\int} e^{i\bar{k}\bar{z}}\Psi_{21}(z,k)\overline{\nu(z)}\;ds(z). 
\end{array}
\end{equation}
All of these computations are performed in practice with $|k|\leq R$ to stabilize the reconstruction in the presence of noise.  The scattering data is set to zero for $|k|>R$.  This approach has been proved to be a nonlinear regularization strategy in the case of real-valued conductivities \cite{Knudsen2009}. 
   Parallel computing can be used to solve equations \eqref{eq:BIEu1u2}-\eqref{eq:scatBIEs} since each of these equations is solved for each $k$ independently.   Further implementation details are found in Section \ref{sec:implementation}.

\vspace{1em}
\noindent\textbf{\underline{Step 2:} Computation of  CGO Solutions:}\\
Let  $\overline{\Omega^+}$ be a domain slightly larger than $\Omega$.  This will be needed to numerically compute the $\dez$ and $\dbarz$ derivatives of the CGO solutions $M(z,0)$ required to form the matrix potential $Q(z)$ in Step 3.
For each $z\in\overline{\Omega^+}$,  solve the $\dbar_k$ equation
\begin{equation}\label{eq:dbarK}
\dbar_k M(z,k)=M(z,\bar{k})\left[\begin{array}{cc}
e(z,\bar{k}) & 0 \\
0 & e(z,-k)
\end{array}\right]S(k),
\end{equation}
using the fundamental solution $\frac{1}{\pi k}$ for the $\dbar_k$ operator, by solving the decoupled systems
\begin{equation}\label{eq:DbarkM-system1}
\small \begin{array}{c}
\left\{\begin{array}{rcl}
M_{11}(z,k)&=& 1 + \frac{1}{\pi k}\ast \left[M_{12}(z,\bar{k})e(z,-k)S_{21}(k)\right]\\
M_{12}(z,k)&=& 0 + \frac{1}{\pi k}\ast \left[M_{11}(z,\bar{k})e(z,\bar{k})S_{12}(k)\right]
\end{array}\right.\\
\end{array}
\end{equation}
\begin{equation}\label{eq:DbarkM-system2}
\small \begin{array}{c}
\left\{\begin{array}{rcl}
M_{21}(z,k)&=& 0 + \frac{1}{\pi k}\ast \left[M_{22}(z,\bar{k})e(z,-k)S_{21}(k)\right]\\
M_{22}(z,k)&=& 1 + \frac{1}{\pi k}\ast \left[M_{21}(z,\bar{k})e(z,\bar{k})S_{12}(k)\right],
\end{array}\right.
\end{array}
\end{equation}
The convolutions $\ast$ take place in $k$ over the disc of radius $R$.

\vspace{1em}
\noindent\textbf{\underline{Step 3:} From CGO Solutions to the Admittivity:}\\
Using the CGO solutions corresponding to $k=0$, compute the potentials (only one is actually needed)
\begin{equation}\label{eq:MtoQ}
\begin{array}{lcl}
Q_{12}(z)&=&\frac{\dez\left[M_{11}(z,0)+M_{12}(z,0)\right]}{M_{22}(z,0)+M_{21}(z,0)}\\
Q_{21}(z)&=&\frac{\dbarz\left[M_{22}(z,0)+M_{21}(z,0)\right]}{M_{11}(z,0)+M_{12}(z,0)}.
\end{array}
\end{equation}
and from these, compute the admittivity $\gamma(z)$ using either
\begin{equation}\label{eq:QtoGamma}
\gamma(z)=\exp\left\{-\frac{2}{\pi \bar{z}}\ast Q_{12}(z)\right\}=\exp\left\{-\frac{2}{\pi {z}}\ast Q_{21}(z)\right\},
\end{equation}
where the convolution in $z$ takes place over $\overline{\Omega}$ since $Q$ has compact support. 

The reader is referred to \cite{Hamilton2013,Hamilton_Thesis_2012,Herrera2015} for developments and implementations of numerical algorithms for complex D-bar EIT imaging.

\subsection{Inclusion of a priori admittivity information} \label{sec:a_priori}
The low pass-filtering (setting $S(k)=0$ for $|k|>R$) in the non-physical scattering domain has an effect similar to that of traditional low-pass filtering in the standard Fourier domain.  As $|k|\to\infty$, the scattering data $S_{12}(k)\approx \frac{i}{\pi}\widehat{Q_{12}}\left(2k_1,2k_2\right)$ and 
$S_{21}(k)\approx  -\frac{i}{\pi}\widehat{Q_{21}}\left(-2k_1,2k_2\right)$, 
and thus, for large $|k|$ the scattering data are essentially Fourier transforms of the potential $Q(z)$.  Hence, it is reasonable to expect a loss of sharp edges in reconstructions of $\gamma(z)$ from the low-pass filtered scattering data.

In practice, the scattering data computed via the boundary integral equations \eqref{eq:scatBIEs} ``blows up'' to $\pm\infty$ as $|k|$ increases, sometimes as early as $|k|=3.5$ in the presence of noise.  Therefore, a natural question arises.  Is it possible to obtain the scattering data $S(k)$ for a larger radius $R_2\geq R$?  While methods based on post-processing D-bar conductivity images have been proposed \cite{Hamilton2014,Hamilton2015}, the work of Alsaker and Mueller \cite{AlsakerMueller2015} is the first D-bar method which directly includes spatial \textit{a priori} information into the nonlinear reconstruction method.  This information is used in the the scattering transform and in the D-bar equation with parameters that can be adjusted to control the amount of influence the prior has on the reconstruction. 

The scattering data is augmented by the scattering data that corresponds to the prior outside the feasible region of computation of the true scattering data.  Denoting the scattering data from the admittivity prior by $\Sprior$, and the feasible region of computation by $|k|\leq R$, we form the new extended scattering data via the formula
\begin{equation}\label{eq:scatNew}
S_{R,R_2}(k):=\begin{cases}
S(k) & |k|\leq R\\
\Sprior(k) & R<|k|\leq R_2\\
0 & R_2<|k|.
\end{cases}
\end{equation}
where $S(k)$ is computed from current and voltage measurements using \eqref{eq:scatBIEs} for $|k|\leq R$.  The truncation radius $R_2$ controls the amount of influence the inclusion of $\Sprior$ has on the reconstruction.  The larger $R_2$, the greater the influence.  When $R_2=R$, there is no inclusion of $\Sprior$.  Note that since $|\Sprior| \rightarrow 0$ as $|k|\rightarrow\infty$, the influence of $\Sprior$ does not grow without bound as $R_2$ increases.

The second place that \textit{a priori} information is included in the reconstruction method is in the integral forms of the D-bar equations, systems \eqref{eq:DbarkM-system1} and \eqref{eq:DbarkM-system2}.
The $+1$ and $+0$ terms in \eqref{eq:DbarkM-system1}, \eqref{eq:DbarkM-system2} arise from terms of the form
\begin{equation}\label{eq:Mint}
\lim_{R\to\infty}\; \frac{1}{\pi R^2}\int_{|k|\leq R} M_{ij}(z,k)\;dk, \quad i,j=1,2,
\end{equation}
whose limits are $0$ for $M_{12}$ and $M_{21}$ and $1$ for $M_{11}$ and $M_{22}$.  Analogously to \cite{AlsakerMueller2015}, to include \textit{a priori} information encoded in the CGO solutions, the terms in \eqref{eq:Mint} are replaced by a weighted integral, which we will denote by 
\begin{equation}\label{eq:DbarkM-replace} 
\Mint\equiv
\begin{array}{c}
\left\{\begin{array}{rcl}
 \alpha + (1-\alpha)\int_{|k|\leq R_2} M^{\mbox{\tiny\textbf{PR}}}_{ij}(z,k)\;dk,  i=j, \\
 0 + (1-\alpha)\int_{|k|\leq R_2} M^{\mbox{\tiny\textbf{PR}}}_{ij}(z,k)\;dk,   i\neq j
\end{array}\right.
\end{array}
\end{equation}
Note, when $\alpha=1$ and $R_2=R$ the method reduces to the original D-bar method of  Subsection~\ref{sec:Dbar} without \textit{a priori} information.


We summarize the steps of the \textit{a priori} method.  The final approximation to the admittivity is denoted by $\gamma_{\New}$.\\



\noindent\textbf{\underline{Step 0:}  Setup:}\\
Compute the DN map $\Lambda_\gamma$ from the voltage and current measurements and determine an admittivity prior $\gprior$.\\

\noindent\textbf{\underline{Step I:} Computation of Scattering Data $S_{R,R_2}$:}\\
Compute the extendend scattering $S_{R,R_2}$ via \eqref{eq:scatNew}.  This involves using Step~1 of Subsection~\ref{sec:Dbar} to compute the traditional scattering data $S(k)$ for $\;|k|\leq R$.  To obtain $\Sprior$ computationally, the smoothed admittivity prior is first used to compute the potential $\Qprior$ via \eqref{eq:Q}.  Then, for $|k|\leq R_2$, the system \eqref{eq:DkM=QM} is solved, and the resulting matrix of CGO solutions is denoted by $\Mprior(\cdot,k)$.  The scattering data $\Sprior(k)$ is then computed via \eqref{eq:ScatDom} using $\Qprior$ and $\Mprior$ in these equations. \\

\noindent\textbf{\underline{Step II:} Computation of  CGO solutions:} \\
Using the extended scattering data $S_{R,R_2}$, solve the systems \eqref{eq:DbarkM-system1} and \eqref{eq:DbarkM-system2} with \eqref{eq:DbarkM-replace} replacing the constant terms $0$ and $1$ to obtain CGO solutions $M_{ij}(z,k),\;i,j=1,2,\; z\in\Omega$, where
\begin{eqnarray*}
M_{11}(z,k)&=& \Minta + \frac{1}{\pi k}\ast \left[M_{12}(z,\bar{k})e(z,-k)S_{21}(k)\right]\\
M_{12}(z,k)&=& \Mintb + \frac{1}{\pi k}\ast \left[M_{11}(z,\bar{k})e(z,\bar{k})S_{12}(k)\right] \\
M_{21}(z,k)&=& \Mintc + \frac{1}{\pi k}\ast \left[M_{22}(z,\bar{k})e(z,-k)S_{21}(k)\right]\\
M_{22}(z,k)&=& \Mintd + \frac{1}{\pi k}\ast \left[M_{21}(z,\bar{k})e(z,\bar{k})S_{12}(k)\right],
\end{eqnarray*}

%

\noindent\textbf{\underline{Step III:} From CGO solutions to the Admittivity $\gamma_{\New}(z)$:}\\
This is computed in the same manner as Step~3 in Subsection~\ref{sec:Dbar} to obtain $\gamma_{\New}(z)$ via  \eqref{eq:MtoQ} using $M_{ij}(z,k),\; i,j=1,2,$ and subsequently \eqref{eq:QtoGamma}.

\section{Simulation and Implementation} 
\subsection{Simulation of Voltage Data}\label{sec:Examples}
The FEM was used to simulate voltages for each of the test problems using the Complete Electrode Model (CEM)  on the chest-shaped domain in Figure \ref{Fig:heartnlung_prior} of perimeter 1016 mm, with $L=32$ electrodes of length $22$ mm and height $13.5$ mm.  The contact impedance was set to $2.4\by 10^{-3}$ on all electrodes, and trigonometric current patterns with amplitude $C=1$mA were used.  The trigonometric current patterns are given by
\begin{equation}\label{eq:trigCPs}
\small T^j_\ell :=\begin{cases}
C\cos\left(j\theta_\ell\right) & 1\leq\ell, \;\; 1\leq j\leq \frac{L}{2}\\
C\sin\left(\left(\frac{L}{2}-j\right)\theta_\ell\right) & 1\leq\ell, \;\; \frac{L}{2} + 1\leq j\leq L-1,
\end{cases}
\end{equation}
where $\theta_\ell=\frac{2\pi\ell}{L}$ corresponds to the angle of the center point $z_\ell=R\left(\theta_\ell\right)e^{i\theta_\ell}$ of the $\ell$-th electrode $e_\ell$.  The quantity $T^j_\ell$ therefore represents the current applied on $e_\ell$ corresponding to the $j$-th current pattern.  Note that $L-1$ linearly independent current patterns were applied since $L$ electrodes were used in the simulations.

Zero mean Gaussian relative noise was added to each complex-valued vector of simulated voltages $V^j$ in the same manner as \cite{Hamilton2013} as follows.  Let $\eta$ denote the the desired noise level and $N^j$ a vector of Gaussian zero mean noise that is unique for each current pattern $j$ (and each test scenario).  Then, the real and imaginary parts of the noisy voltage data $\widetilde{V}^j$ were computed as 
\begin{equation}\label{eq:noisyVolts}
\begin{array}{lcr}
\Re\left(\widetilde{V}^j\right)&=& \Re\left({V}^j\right)+ \eta\; \max\left|\Re\left({V}^j\right)\right|\;N^j\\
\Im\left(\widetilde{V}^j\right)&=& \Im\left({V}^j\right)+ \eta\; \max\left|\Im\left({V}^j\right)\right|\;N^j.
\end{array}
\end{equation}
The discrete approximation $\Lambda^M_\gamma$ to the D-N map was computed as in \cite{Isaacson2004,Isaacson2006},  which we summarize briefly here.  Denoting by $t^j_{\ell}$ the $({\ell},j)$-th entry of the matrix of applied currents with each column normalized with respect to the $l^2$-vector norm, $t^j_{\ell} = \frac{T^j}{\|T^j\|_2}$, let $v^j_{\ell}$ denote the entries of the $j$-th voltage vector normalized so that $\sum_{\ell=1}^Lv_{\ell}^j = 0$ and $v_{\ell}^j = \frac{V_{\ell}}{\|T^j\|_2}$.  Let $\left|e_\ell\right|$ denote the area of the $\ell$-th electrode.  Then $\Lambda^M_\gamma  =(R^M_\gamma)^{-1}$ where the $(m,n)$-th entry of $R_\gamma^M$ is given by
\begin{equation} \label{Rmn}
R_\gamma^M(m,n) := \frac{\gamma_0}{\left|e_\ell\right|}\sum_{\ell=1}^L t_{\ell}^m v_{\ell}^n.
\end{equation}

\begin{figure}[!h]
\centering
{\includegraphics[width=100pt]{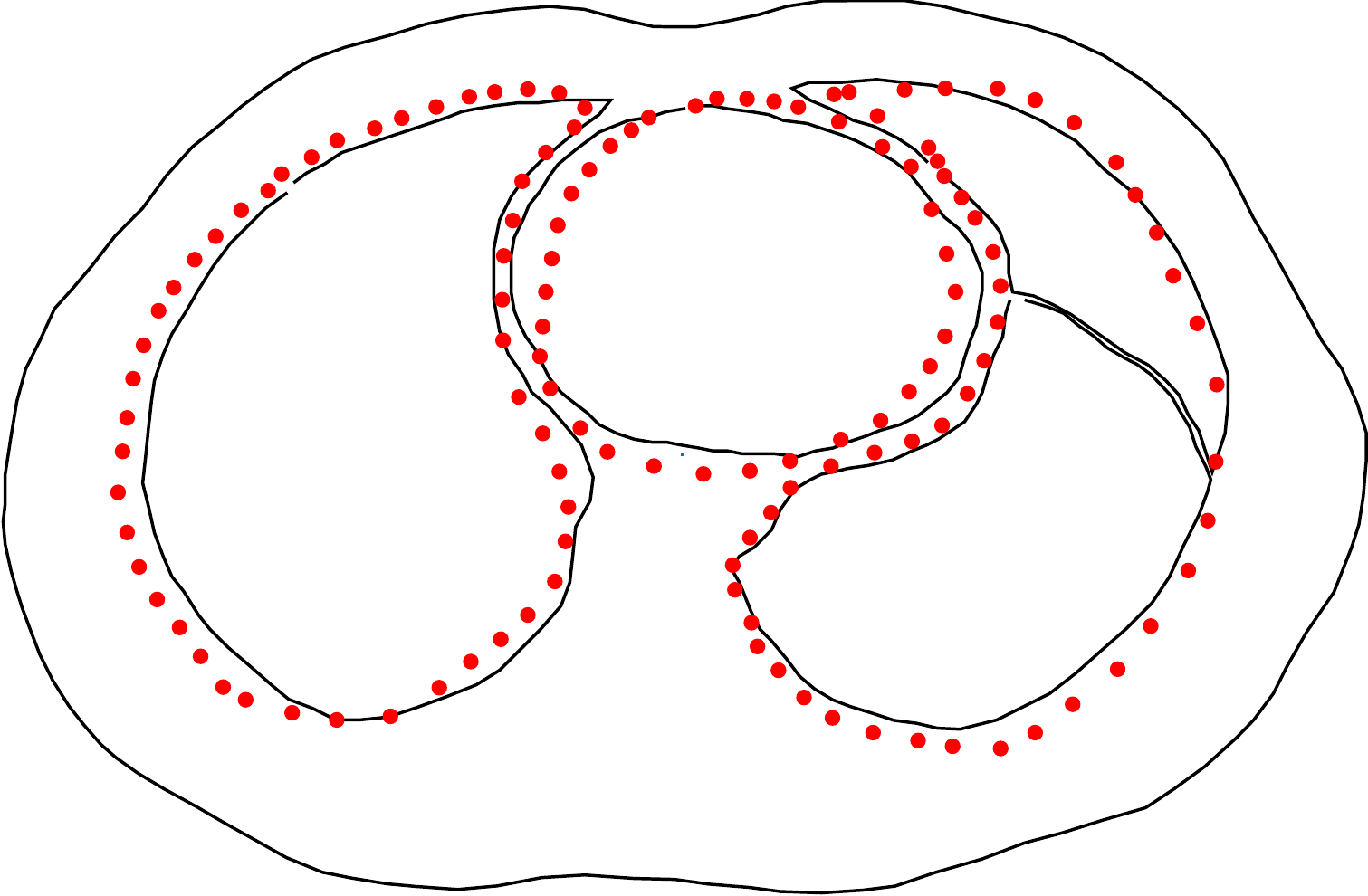}}\quad
{\includegraphics[width=100pt]{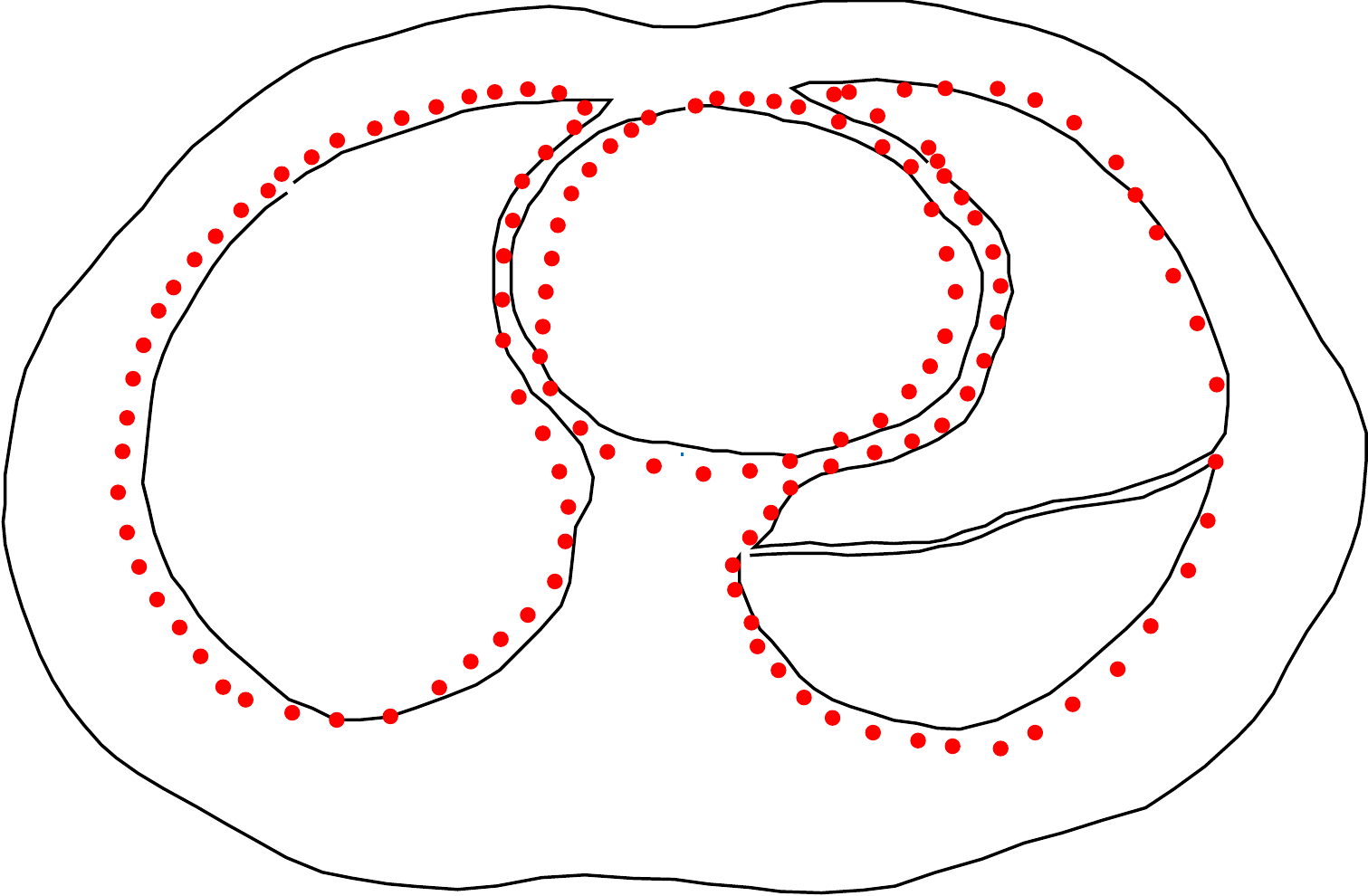}}\\
(a)\hspace*{100pt} (b)
\caption{\label{fig:phantoms} The test examples studied simulate two pathologies:  (a) a pneumothorax in the ventral part of the left lung and (b) a pleural effusion in the dorsal part of the left lung.  The black lines correspond to true boundaries in the simulations, and the superimposed red dots are the organ boundaries used in the construction of the admittivity prior $\gprior$ before smoothing.}
\end{figure}

Figure~\ref{fig:phantoms} shows the two simulations: (a) a pneumothorax, (b) a pleural effusion.  For both sets of simulated data, the admittivity of the heart was $1.1+0.6i$ S/m, the lungs $0.5+0.4i$ S/m, and the background $0.8+0.4i$ S/m.  The pneumothorax was set to $0.25+0i$ S/m and the pleural effusion to $1.1+0.6i$ S/m.

\subsection{Implementation of the a priori method} \label{sec:implementation}

In this paper, the admittivity prior $\gprior$ was computed using a standard D-bar reconstruction $\gamma_{\mbox{\tiny{{DB}}}}$ recovered using Steps~1-3 of Section~\ref{sec:Dbar} with the measured data $\Lambda_\gamma$.  However, in practice, any initial prior $\gprior$ can be used, making the method easily adaptable to other approaches.  \\

\noindent\textbf{\underline{Step 0:}}  The matrix approximation to the DN map $\Lambda_\gamma$ was formed using the noisy voltages computed from the CEM.   The admittivity prior $\gprior$ was formed as follows.  First the standard D-bar reconstruction $\gDB$ was computed using Steps~1-3 of Section~\ref{sec:Dbar} (see \cite{Hamilton2013} for details regarding the computation of $\gDB$).  Next, using the spatial heart and lungs prior (see the red dots of Figure~\ref{Fig:heartnlung_prior}), the average value of the pixels in each region (heart, left lung, right lung, and background) were computed and the corresponding average assigned to each region to form the admittivity prior $\gprior$.  Note that the spatial prior does not assume any pathology is present.   The prior $\gprior$ was then mollified to a $C^1$ smooth version and $\Qprior$ computed using finite differences for the $\dez$ and $\dbarz$ derivatives of $\log\left(\gprior(z)\right)$.

\noindent\textbf{\underline{Step I:}}  The extended scattering data $S_{R,R_2}$ was computed via \eqref{eq:scatNew}.  Using the DN maps $\Lambda_\gamma$ and $\Lambda_1$, the traditional scattering data $S(k)$ for $|k|\leq R$ was determined via Step~1 of Section~\ref{sec:Dbar}.   The reader is referred to \cite{Hamilton2013} for the computational details of computing $u_1$ and $u_2$ and subsequently $\psi_{12}$ and $\psi_{21}$.  Briefly, 
the Fredholm integral equations for $u_1$ and $u_2$ \eqref{eq:BIEu1u2} are solved by a Galerkin method, and the integrals for evaluating $\psi_{12}$ and $\psi_{21}$ and scattering data $S(k)$, $|k|\leq R$  in \eqref{eq:scatBIEs}  are computed using a Simpson's rule.  The scattering prior $\Sprior$ is determined as follows.  First, the admittivity prior $\gprior$ is smoothed to compute the potential $\Qprior$ via \eqref{eq:Q}.  Then, for $|k|\leq R_2$, the system \eqref{eq:DkM=QM} is solved for $M^{\prior}$ using Fourier transforms on the following two decoupled systems:
\begin{equation}\label{eq:DkM-QM-systems-Prior}
\begin{array}{c}
\left\{\begin{array}{rcl}
M^{\prior}_{11}(z,k)&=& 1 + \frac{1}{\pi z} \ast \left[Q^{\prior}_{12}(z)M^{\prior}_{21}(z,k)\right]\\

M^{\prior}_{21}(z,k)&=& 0 + \frac{e(z,-k)}{\pi\bar{z}}\ast\left[Q^{\prior}_{21}(z)M^{\prior}_{11}(z,k)\right]
\end{array}\right.\\
\\
\left\{\begin{array}{rcl}
M^{\prior}_{12}(z,k)&=& 0 + \frac{e(z,\bar{k})}{\pi z}\ast\left[Q^{\prior}_{12}(z)M^{\prior}_{22}(z,k)\right]\\

M^{\prior}_{22}(z,k)&=& 1 + \frac{1}{\pi \bar{z}}\ast\left[Q^{\prior}_{21}(z)M^{\prior}_{12}(z,k)\right]
\end{array}\right.,
\end{array}
\end{equation}
where the convolutions take place in $z$ over ${\Omega}$.  Using a uniform $z$-grid of size $2^m\by 2^m$ with stepsize $h$, convolutions such as $\frac{1}{\pi z}\ast f(z)$ can be implemented as
\[\footnotesize \frac{1}{\pi z}\ast f(z)=h^2\texttt{IFFT2}\left(\texttt{FFT2}\left(\frac{1}{\pi z}\right)\texttt{FFT2}\left(f(z)\right)\right).\]
The scattering prior is then evaluated via \eqref{eq:ScatDom} using Simpson's rule, and the combined scattering $S_{R,R_2}$ is subsequently formed via \eqref{eq:scatNew}.\\

\noindent\textbf{\underline{Step II:}}
Choose a regularization weight $\alpha\in[0,1]$.  Using the combined scattering data $S_{R,R_2}$, the CGO solutions $M_{11}^{R_2,\alpha}$ and $M_{12}^{R_2,\alpha}$ were recovered using Fourier transforms to solve the modified equations \eqref{eq:DbarkM-systems-Prior-Regularized} 
\begin{equation}\label{eq:DbarkM-systems-Prior-Regularized}
\footnotesize
\begin{array}{c}
\left\{\begin{array}{lcc}
M_{11}^{R_2,\alpha}(z,k)= M^{\mbox{\tiny{\textbf{int}}}}_{11}(z)
+ \frac{1}{\pi k}\ast \left[M_{12}^{R_2,\alpha}(z,\bar{k})e(z,-k)S^{21}_{R,R_2}(k)
\right]
\\
&&\\

M_{12}^{R_2,\alpha}(z,k)= M^{\mbox{\tiny{\textbf{int}}}}_{12}(z) +\frac{1}{\pi k}\ast \left[M_{11}^{R_2,\alpha}(z,\bar{k})e(z,\bar{k})S^{12}_{R,R_2}(k)
\right]\end{array}\right.
\end{array}
\end{equation}
where the convolutions take place in $k$ over $|k|\leq R_2$ and  $M^{\mbox{\tiny{\textbf{int}}}}_{ij}$ is computed from \eqref{eq:DbarkM-replace} using a Simpson's rule.  An analogous system is solved to recover $M_{21}^{R_2,\alpha}$ and $M_{22}^{R_2,\alpha}$. \\ 

\noindent\textbf{\underline{Step III:} } The new admittivity is recovered in the same manner as Step~3 of Subsection~\ref{sec:Dbar} to obtain $\gamma_{\New}(z)$ via  \eqref{eq:MtoQ} using finite differences on $M_{ij}^{R_2,\alpha}(z,0),\; i,j=1,2,$ and subsequently Fourier transforms to solve \eqref{eq:QtoGamma}.

\subsection{Examples}
In this work, two noise levels were considered: $0.1\%$ added relative noise and $1.0\%$ relative noise.  For each example, we present results with three values of the truncation radius $R_2$ in the prior, and three regularization weights for the D-bar equation: $\alpha=0,0.5,1$.  Recall that $\alpha=0$ corresponds to the strongest weight and $\alpha=1$ to no weight given (see \eqref{eq:DbarkM-replace}).   Due to the ill-posedness of the inverse problem, the radii $R$ of admissible scattering data is problem specific, and the scattering transform will blow up in the presence of noise at a rate that is more rapid in some directions in the $k$-plane than others.  The value chosen for each example was chosen empirically to be as large as possible without exhibiting blow up in the initial reconstruction without \textit{a priori} information.  The blow-up was more rapid in the case of 1\% noise, and so in those examples a non-uniform truncation of the scattering transform was used.  In such cases a threshold of the scattering data $S(k)$ was enforced by setting $S_{ij}(k)=0$ if $\Re(\left|S_{ij}\right|)>0.15 $ or $\Im(\left|S_{ij}\right|)>0.15$, where the value $0.15$ was chosen empirically to be the largest permissible value of the magnitude.  Determining such a threshold is intuitive from a plot of the scattering data since the blowup rate is exponential.


The admittivity prior $\gprior$ consisted of approximate knowledge of the organ boundaries (see Figure~\ref{Fig:heartnlung_prior}) with no assumption of pathology in the lungs.     These average values for the prior are given in  Tables \ref{table:pneumo} and \ref{table:fluid}.


\begin{table*}[t] 
\caption{Admittivity values used in the example of a simulated pneumothorax in the left lung.} 
  \centering
  \begin{tabular}{*{6}{c}}
\hline
Admittivities & Background& Left Lung & Right Lung & Heart &Pneumothorax \\
\hline \hline
Truth  & $0.8+0.4i$  & $0.5+0.2i$ & $0.5+0.2i$ & $1.1+0.6i$ & $0.25+0i$ \\
Prior 0.1\% Noise &     0.79 + 0.40i &   0.66 + 0.28i &  0.64 + 0.29i &  0.89 + 0.48i  & N/A \\
Prior 1.0\% Noise & 0.79+0.39i & 0.66 + 0.25i  & 0.64 + 0.28i&  0.84 + 0.47i & N/A \\
\hline
    \hline
  \end{tabular}
  \label{table:pneumo}
\end{table*}

\begin{table*}[t] 
\caption{Admittivity values used in the example of a simulated pleural effusion in the left lung.} 
  \centering
  \begin{tabular}{*{6}{c}}
\hline
Admittivities & Background& Left Lung & Right Lung & Heart &Pleural Effusion \\\hline \hline
Truth  & $0.8+0.4i$  & $0.5+0.2i$ & $0.5+0.2i$ & $1.1+0.6i$ & $1.1+0.6i$ \\
Prior 0.1\% Noise & 0.80+0.40i &   0.77 + 0.39i  & 0.62 + 0.29i  & 0.92 + 0.47i  & N/A \\
Prior 1.0\% Noise &  0.79 + 0.40i &  0.74 + 0.39i  &  0.59 + 0.31i   &  0.91 + 0.52i & N/A \\
\hline
    \hline
  \end{tabular}
  \label{table:fluid}
\end{table*}

\subsubsection{Example 1:  Simulated Pneumothorax:}
\begin{figure}[!h]
\begin{center}
 {\bf Conductivity}\\
 {\tiny$\mathbf{\alpha=1.00}$ \hspace{0.85in} \tiny$\mathbf{\alpha=0.50}$ \hspace{0.85in} \tiny$\mathbf{\alpha=0.00}$}\\
 
\vspace{.25em}
{\includegraphics[width=300pt]{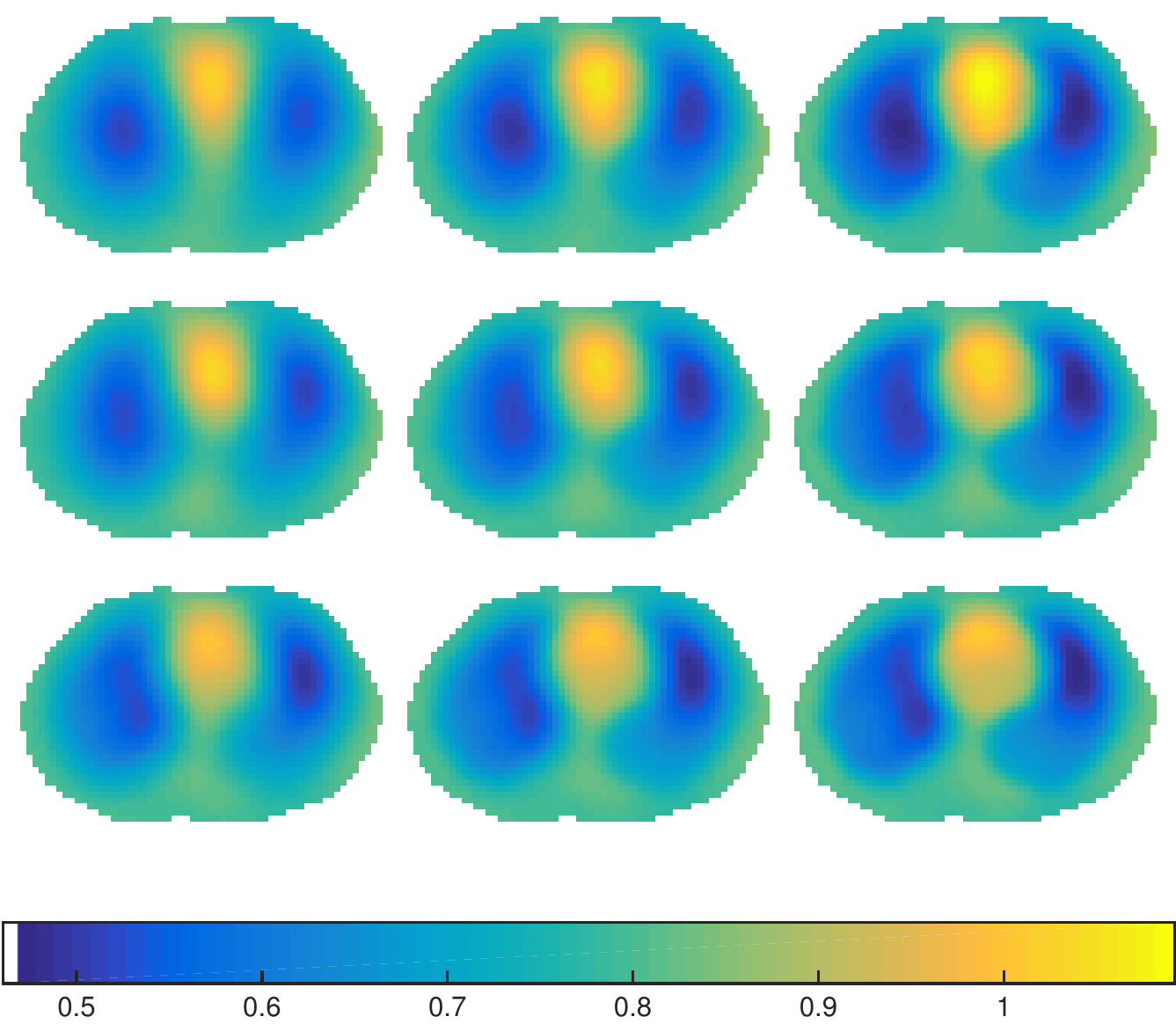}}
\end{center}
\begin{center}
{\bf Permittivity}\\
 {\tiny$\mathbf{\alpha=1.00}$ \hspace{0.85in} \tiny$\mathbf{\alpha=0.50}$ \hspace{0.85in} \tiny$\mathbf{\alpha=0.00}$}\\
 
\vspace{.25em}
{\includegraphics[width=300pt]{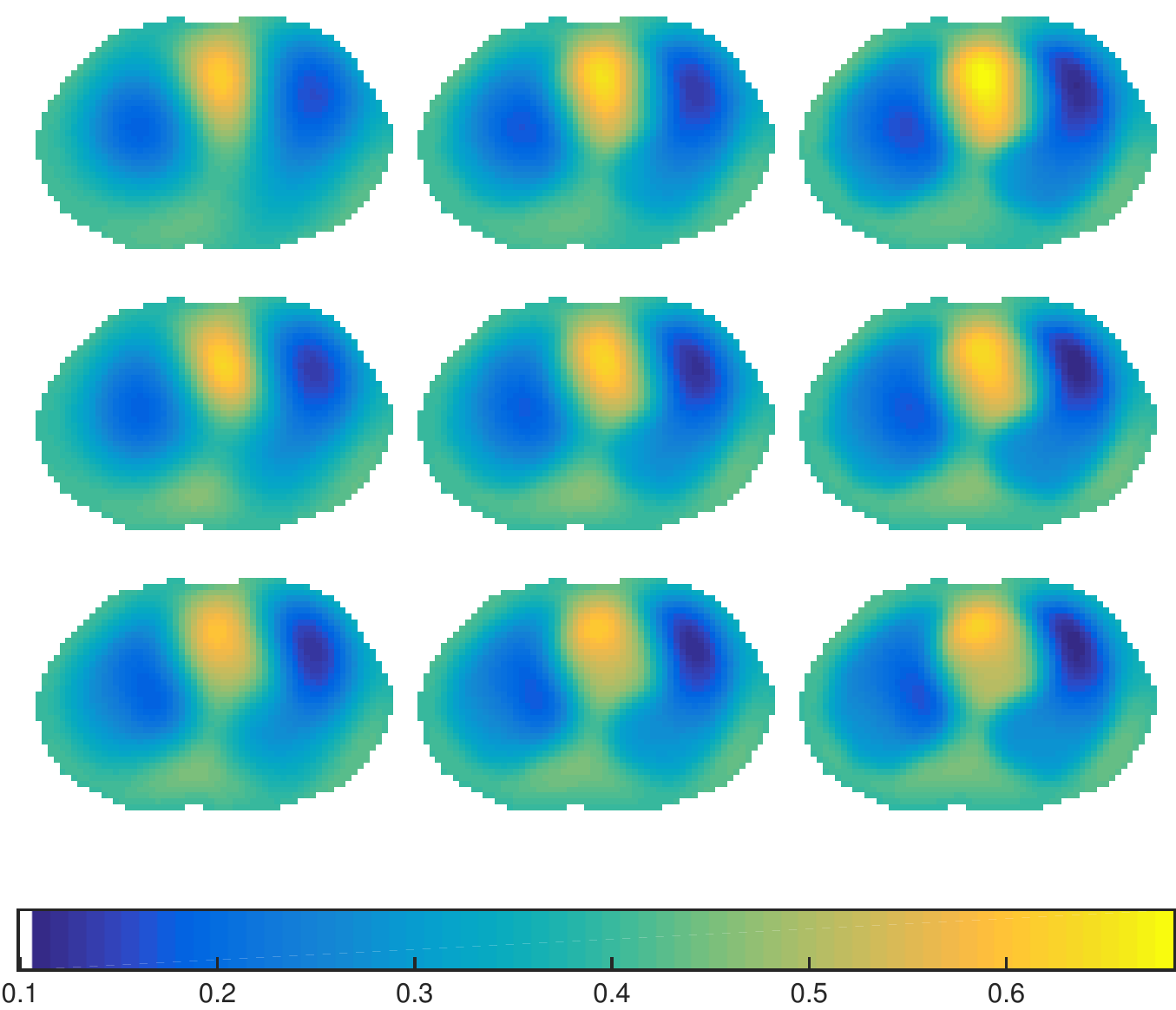}}
\end{center}
\vspace{-1em}
\caption{\label{fig:pneumo-0o1Noise} Reconstructions of simulated pneumothorax with 0.1\% added noise. Regularization parameter  $\alpha=1, 0.5,0$ increases the influence of $\Mint$ as $\alpha$ decreases, and  $R_2=4.5,6.5,8.5$ (rows) increases the influence of $\Sprior$ as $R_2$ increases. No pneumothorax is assumed to be present in the prior.}
\end{figure}

\begin{figure}[!h]
\begin{center}
 {\bf Conductivity}\\
 {\tiny$\mathbf{\alpha=1.00}$ \hspace{0.85in} \tiny$\mathbf{\alpha=0.50}$ \hspace{0.85in}\tiny$\mathbf{\alpha=0.00}$}\\
 
\vspace{.25em}
{\includegraphics[width=300pt]{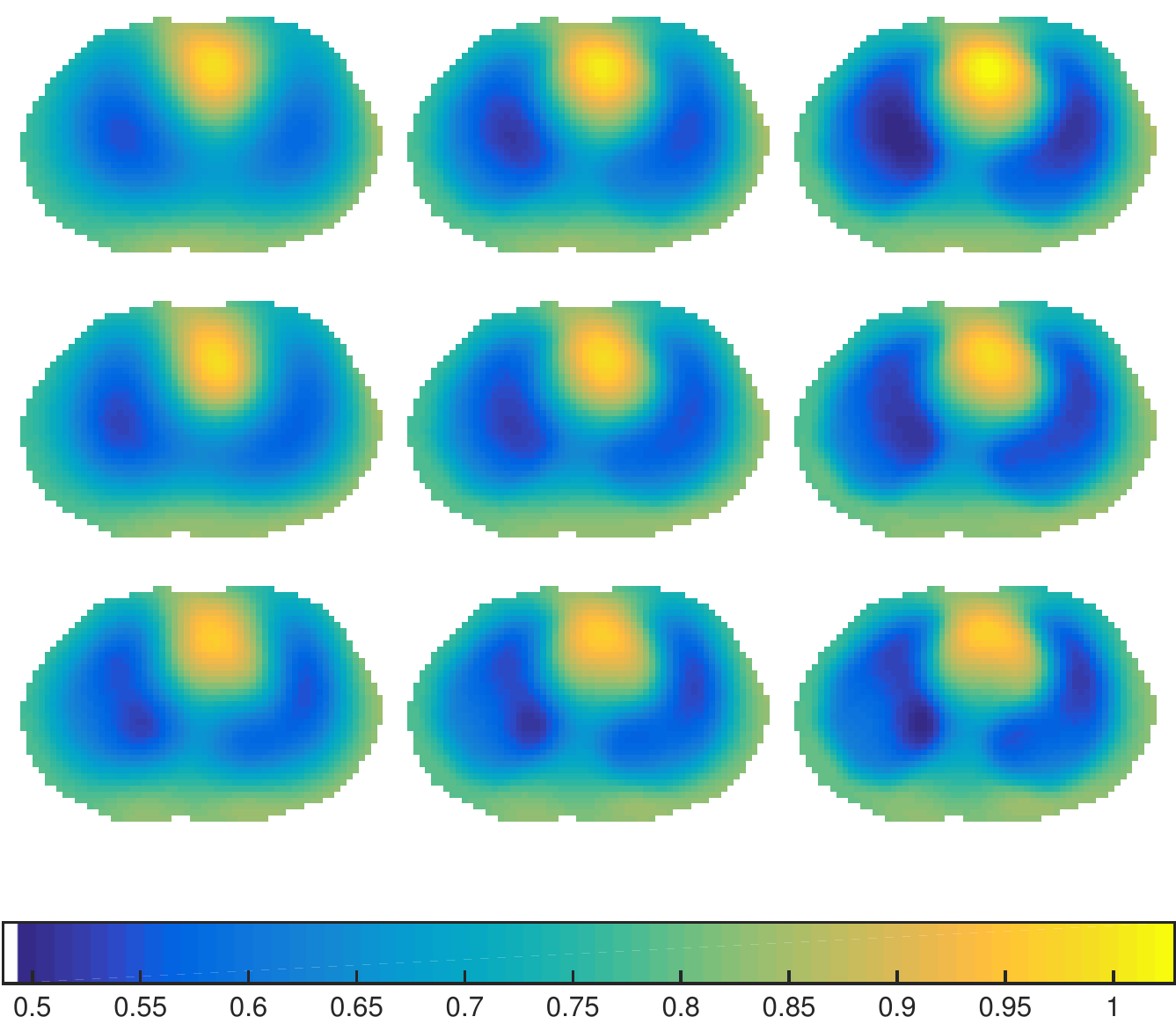}}
\end{center}
\begin{center}
{\bf Permittivity}\\
 {\tiny$\mathbf{\alpha=1.00}$ \hspace{0.85in} \tiny$\mathbf{\alpha=0.50}$ \hspace{0.85in} \tiny$\mathbf{\alpha=0.00}$}\\
 
\vspace{.25em}
{\includegraphics[width=300pt]{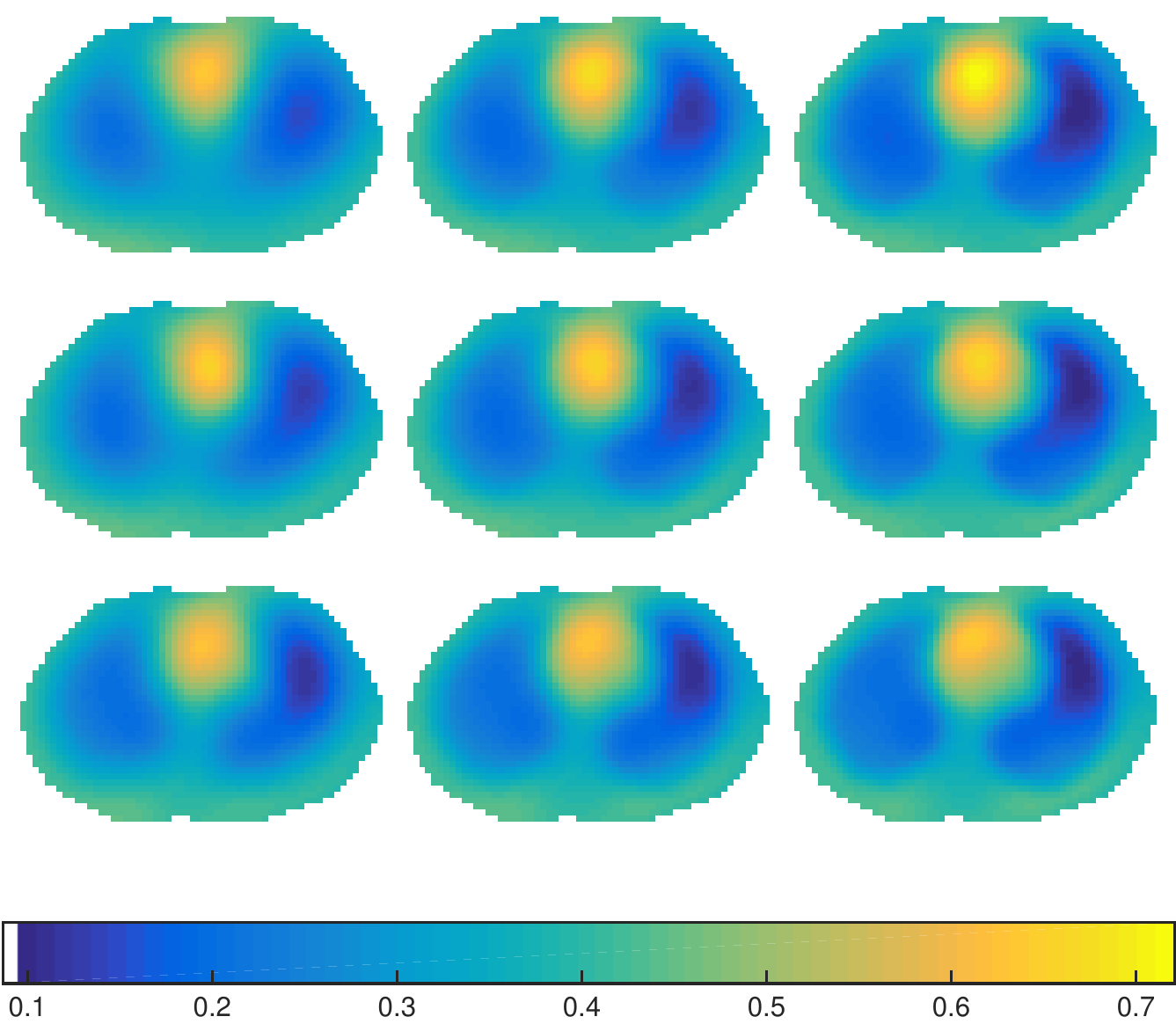}}
\end{center}
\vspace{-1em}
\caption{\label{fig:pneumo-1o0Noise} Reconstructions of simulated pneumothorax with 1.0\% added noise. Regularization parameter  $\alpha=1, 0.5,0$ increases the influence of $\Mint$ as $\alpha$ decreases, and  $R_2=4,6,8$ (rows) increases the influence of $\Sprior$ as $R_2$ increases. No pneumothorax is assumed to be present in the prior.}
\end{figure}

This test problem corresponds to phantom (a) in Figure~\ref{fig:phantoms}.  The preliminary reconstruction with no prior was computed for the 0.1\% added noise case using a radius of $R=4.5$,  and for the 1\% added noise case using a nonuniform truncation with a maximum radius of $R=4.0$.   Table~\ref{table:pneumo} contains the values of the true admittivity in each region as well as the values assigned to the heart and lung prior for 0.1\% and 1\% noise.  We emphasize that we assume only approximate knowledge of the boundaries of the heart and lungs (see the red dots in Figure~\ref{fig:phantoms} (a)), and no knowledge of the presence of a pneumothorax.    Reconstructions for the 0.1\% added noise case with truncation radii for the prior $R_2=4.5, 6.5, 8.5$ and weights $\alpha= 0, 0.5, 1$ are found in Figure \ref{fig:pneumo-0o1Noise}.  Reconstructions for the $1\%$ added relative noise case with truncation radii for the prior $R_2=4, 6, 8$ and weights $\alpha= 0, 0.5, 1$ are shown in  Figure~\ref{fig:pneumo-1o0Noise}.

\subsubsection{Example 2:  Simulated Pleural Effusion }
\begin{figure}[!h]
\begin{center}
 {\bf Conductivity}\\
 {\tiny$\mathbf{\alpha=1.00}$ \hspace{0.85in} \tiny$\mathbf{\alpha=0.50}$ \hspace{0.85in} \tiny$\mathbf{\alpha=0.00}$}\\
 
\vspace{.25em}
{\includegraphics[width=300pt]{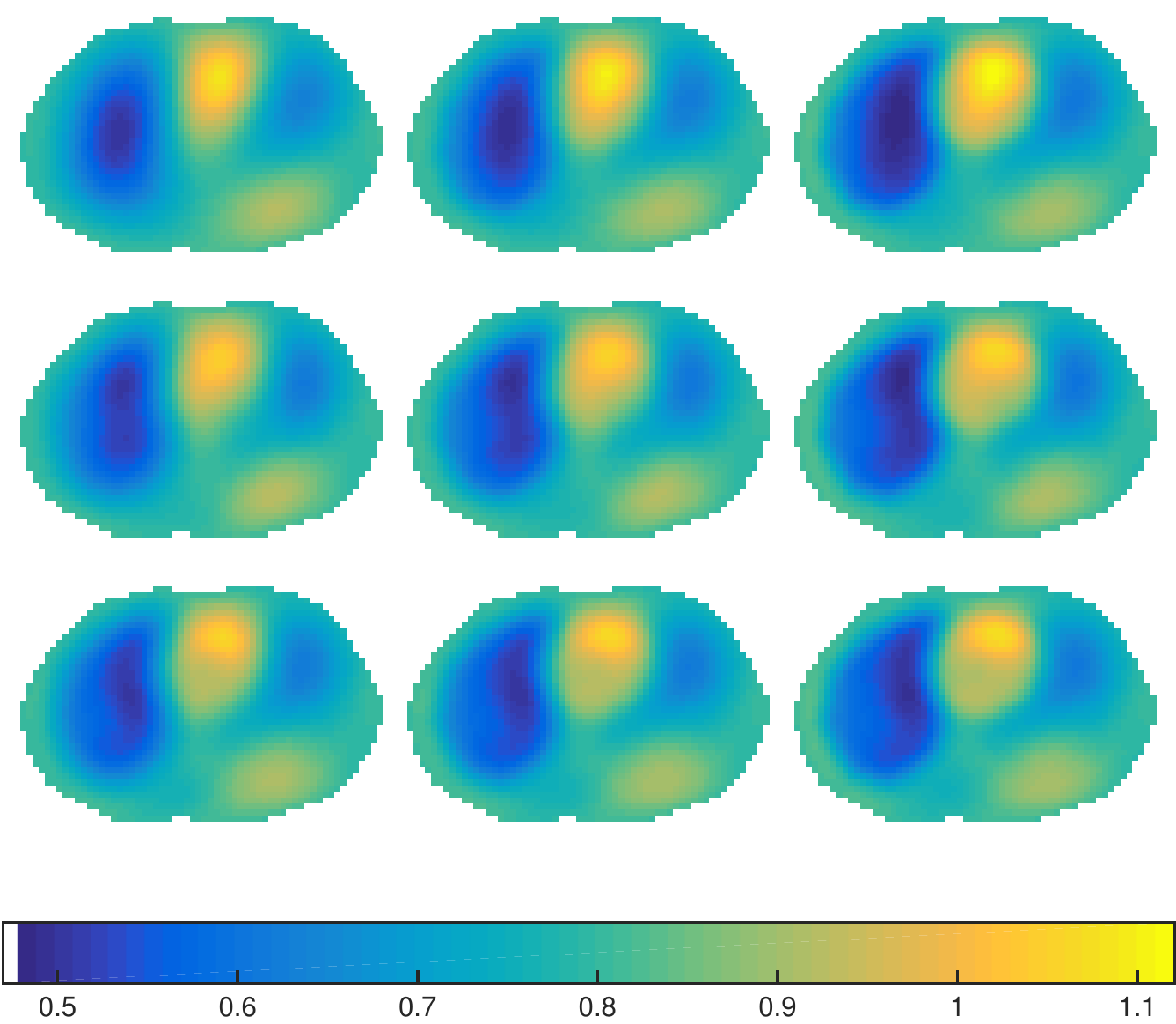}}
\end{center}
\begin{center}
{\bf Permittivity}\\
 {\tiny$\mathbf{\alpha=1.00}$ \hspace{0.85in} \tiny$\mathbf{\alpha=0.50}$ \hspace{0.85in} \tiny$\mathbf{\alpha=0.00}$}\\
 
\vspace{.25em}
{\includegraphics[width=300pt]{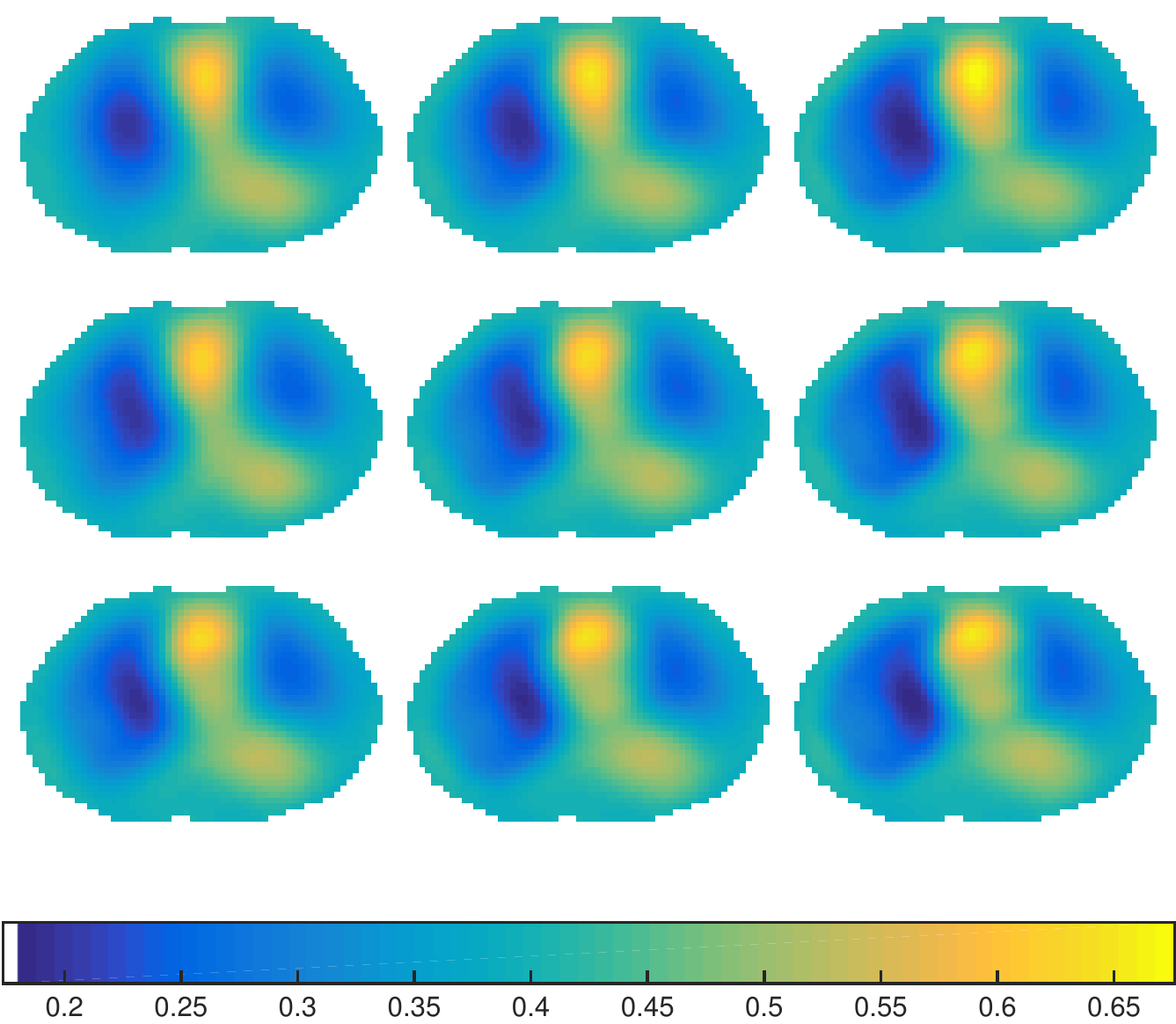}}
\end{center}
\vspace{-1em}
\caption{\label{fig:fluid-0o1Noise} Reconstructions of simulated pleural effusion with 0.1\% added noise. Regularization parameter $\alpha=1, 0.5,0$ increases the influence of $\Mint$ as $\alpha$ decreases, and  $R_2=5.5,8,11$ (rows) increases the influence of $\Sprior$ as $R_2$ increases. No effusion is assumed to be present in the prior.}
\end{figure}

\begin{figure}[!h]
\begin{center}
 {\bf Conductivity}\\
 {\tiny$\mathbf{\alpha=1.00}$ \hspace{0.85in} \tiny$\mathbf{\alpha=0.50}$ \hspace{0.85in} \tiny$\mathbf{\alpha=0.00}$}\\

\vspace{.25em}
{\includegraphics[width=300pt]{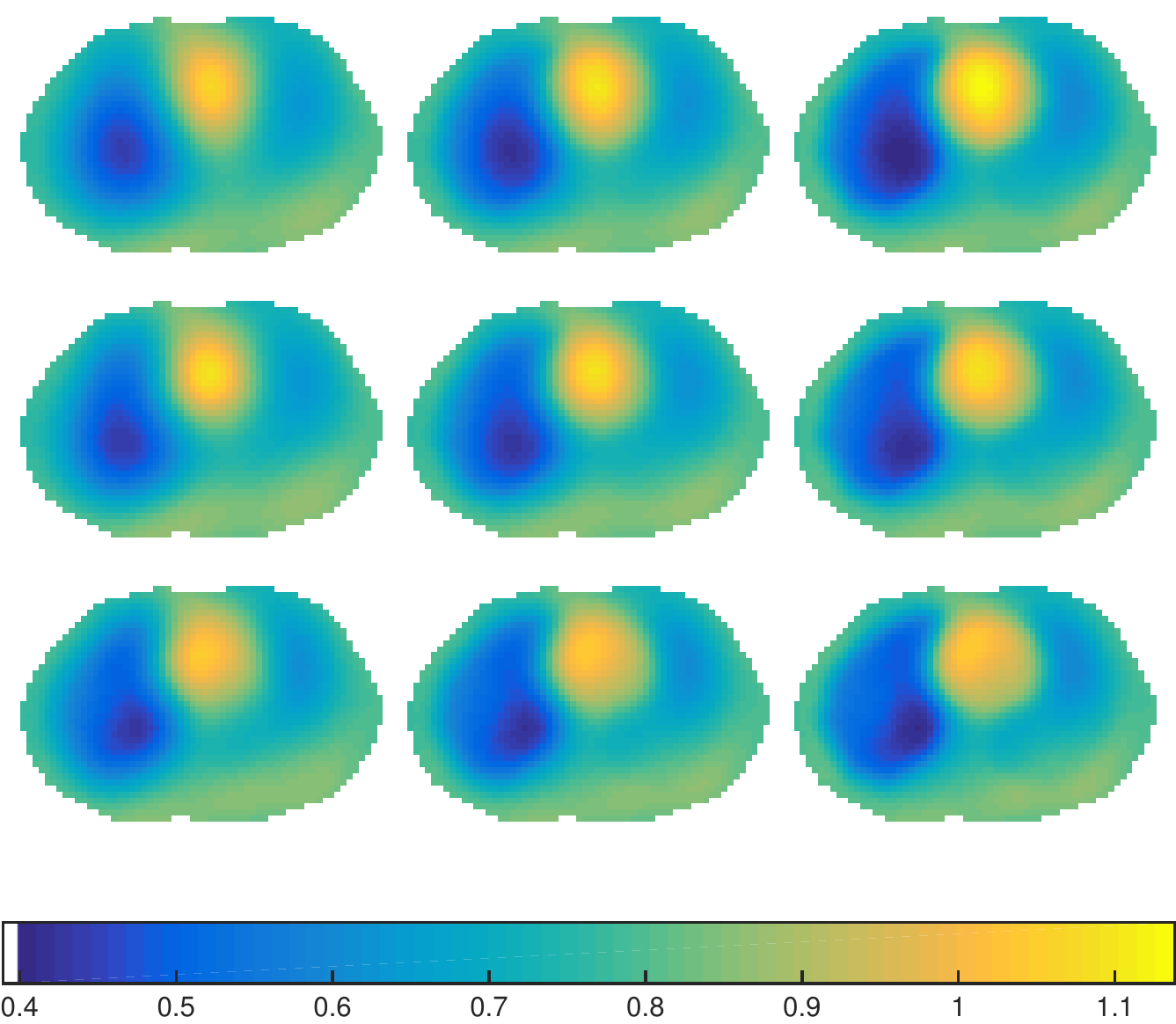}} 
\end{center}
\begin{center}
{\bf Permittivity}\\
 {\tiny$\mathbf{\alpha=1.00}$ \hspace{0.85in} \tiny$\mathbf{\alpha=0.50}$ \hspace{0.85in} \tiny$\mathbf{\alpha=0.00}$}\\

\vspace{.25em}
{\includegraphics[width=300pt]{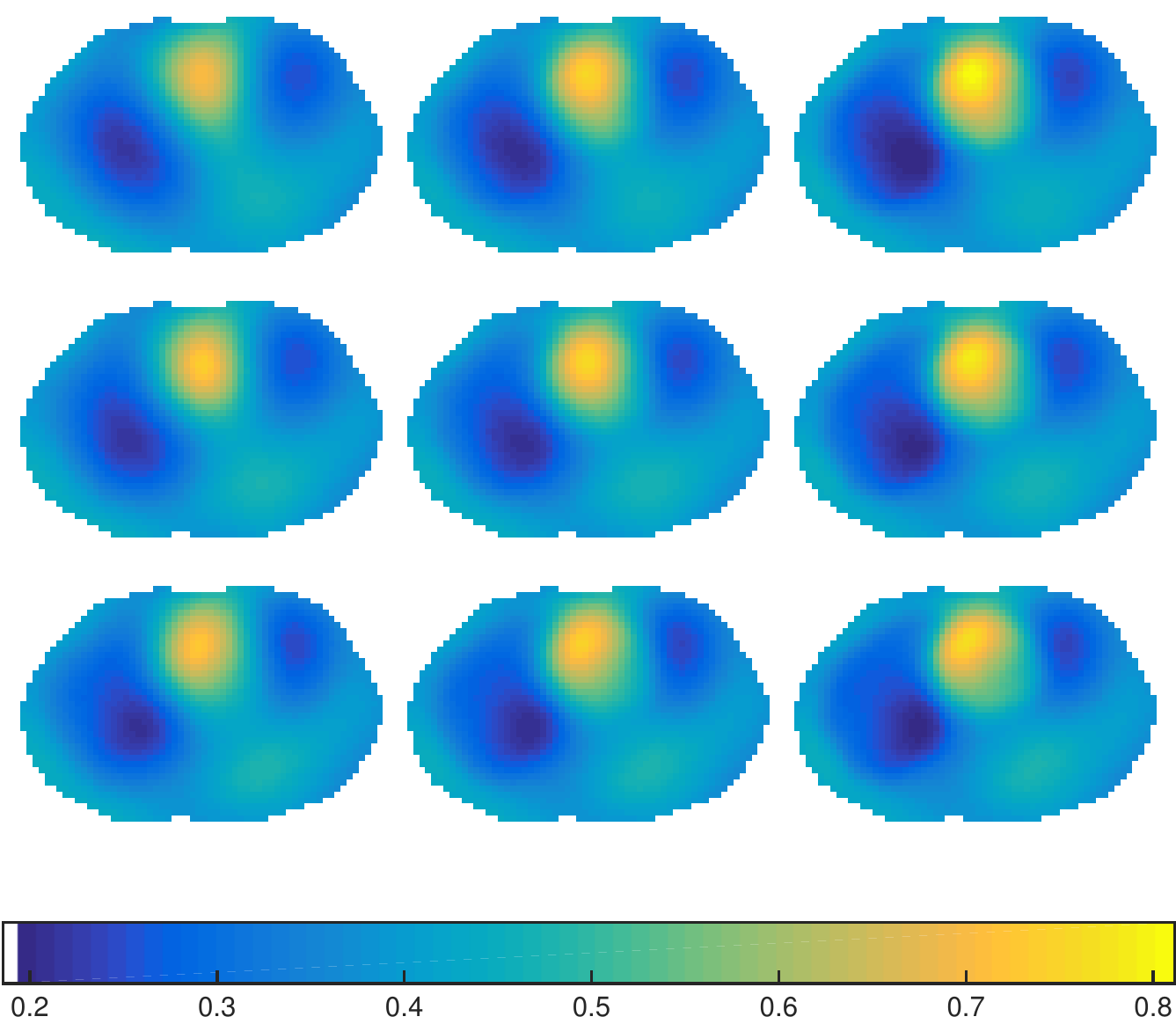}}
\end{center}
\vspace{-1em}
\caption{\label{fig:fluid-1o0Noise} Reconstructions of simulated pleural effusion with 1.0\% added noise. Regularization parameter  $\alpha=1, 0.5,0$ increases the influence of $\Mint$ as $\alpha$ decreases, and  $R_2=4.5,6.5,8.5$ (rows) increases the influence of $\Sprior$ as $R_2$ increases. No effusion is assumed to be present in the prior.}
\end{figure}

This test problem corresponds to phantom (b) in Figure~\ref{fig:phantoms}.  The preliminary reconstruction with no prior was computed for the 0.1\% added noise case using a radius of  $R=5.5$ and for the 1\% added noise case using a nonuniform truncation with a maximum radius of $R=4.5$.  Table~\ref{table:fluid} presents the average values used in the prior $\gprior$ for each noise level.  Reconstructions for the 0.1\% added noise case with truncation radii for the prior $R_2=5.5, 8, 11$ and weight $\alpha= 0, 0.5, 1$ are found in Figure~\ref{fig:fluid-0o1Noise}.  Reconstructions for the $1\%$ added relative noise case with truncation radii for the prior $R_2=4.5, 6.5, 8.5$ and weight $\alpha= 0, 0.5, 1$ are shown in  Figure~\ref{fig:fluid-1o0Noise}.


\clearpage

\section{Discussion and Conclusions}\label{Sec:discussion}
The reader is advised to view the images on a computer screen if possible, since details in the color map are likely masked in printed versions.

In Figures \ref{fig:pneumo-0o1Noise}, \ref{fig:pneumo-1o0Noise}, \ref{fig:fluid-0o1Noise}, and \ref{fig:fluid-1o0Noise}, the upper left figure is the same as the preliminary reconstruction (ie, no prior), and it is evident that the spatial resolution of the organ boundaries improves with the introduction of the prior and as the influence of the prior increases.  
In the case of the pneumothorax, no pathology is evident in the preliminary reconstruction, but as the influence of the prior increases, even though the prior includes no assumption of pathology, the pneumothorax is clearly visible in the reconstructions.  However, in both the conductivity and permittivity images, a lower conductivity and permittivity region becomes evident in the dorsal right lung as well, which is an artifact of the reconstruction, and it becomes stronger as the weighting of $\Mint$ increases ($\alpha = 0.5$ and $0$.)  This artifact is less pronounced in the permittivity images,  and is arguably not present in the $1.0\%$ added noise case in the permittivity images.

The presence of the simulated pleural effusion, on the other hand, is clearly evident in the preliminary reconstructions  for both conductivity and permittivity and for both noise levels.  The presence of the prior improves the spatial resolution of the organs and the actual conductivity and permittivity values in the region of the effusion, but since the regularization results in reconstructed conductivity and permittivity functions that are smooth, there is a smooth transition from the healthy ventral portion of the left lung to the effusion, and so the boundary is far from as sharp as in the piecewise constant phantom.  In practice, image segmentation is often used on reconstructed EIT images, which would likely improve the appearance of the reconstructed images.  Alternatively, once a pathology is visible, an iterative method could then be invoked as in \cite{AlsakerMueller2015} which segments the prior in the region of a possible pathology potentially sharpening the pathology even more.  Post-processing approaches are left for future work.

\begin{figure}[!h] 

\centering
\rotatebox{90}{\textbf{Conductivity}}
{\includegraphics[width=100pt]{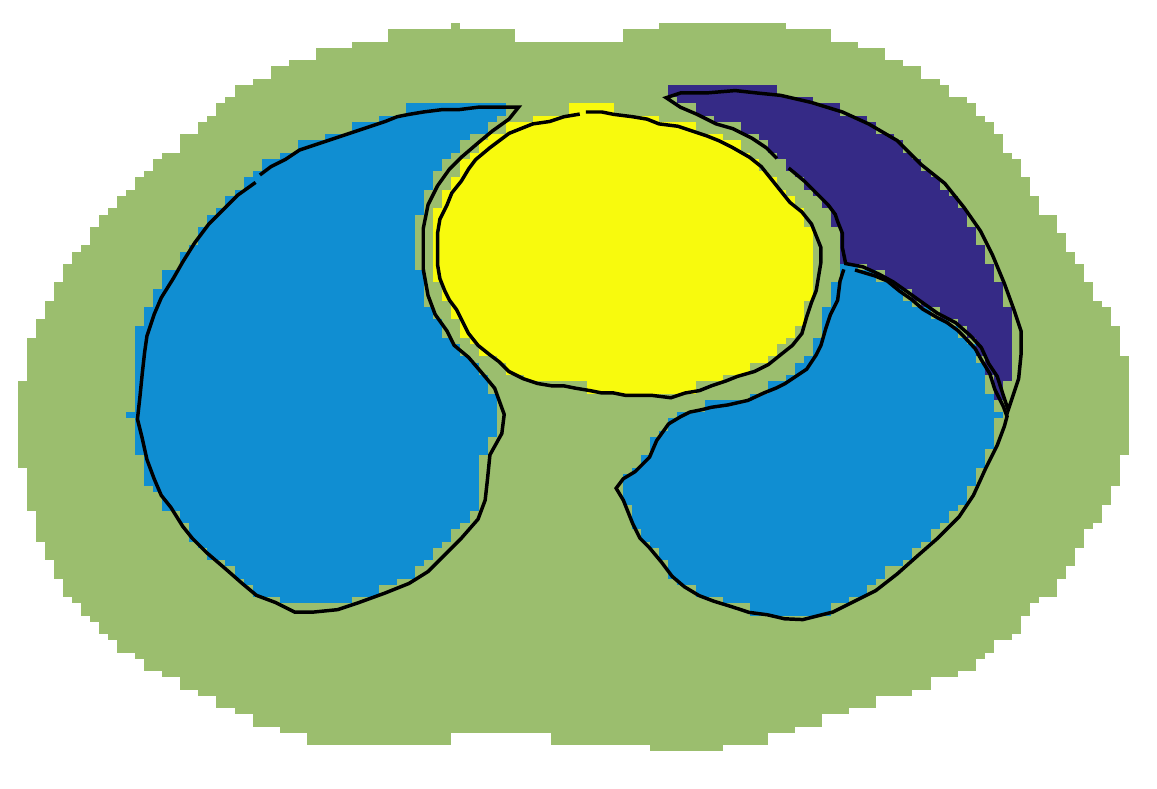}}
{\includegraphics[width=100pt]{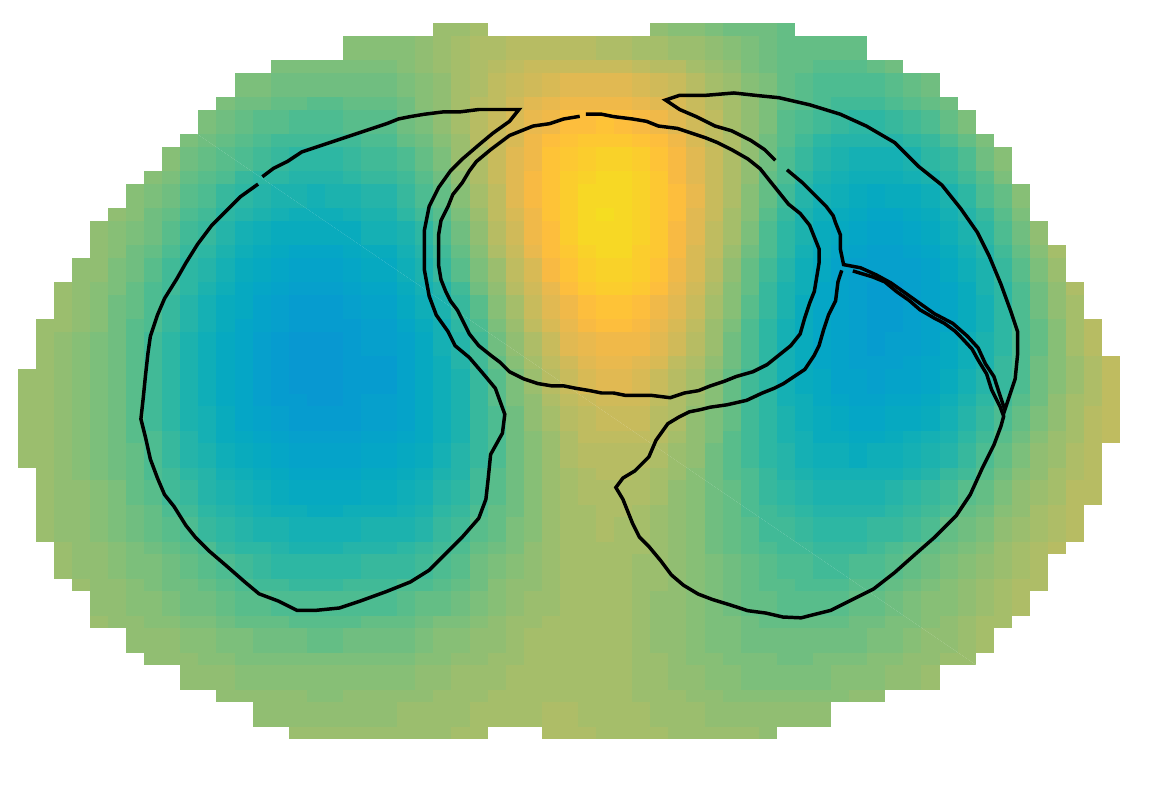}}
{\includegraphics[width=100pt]{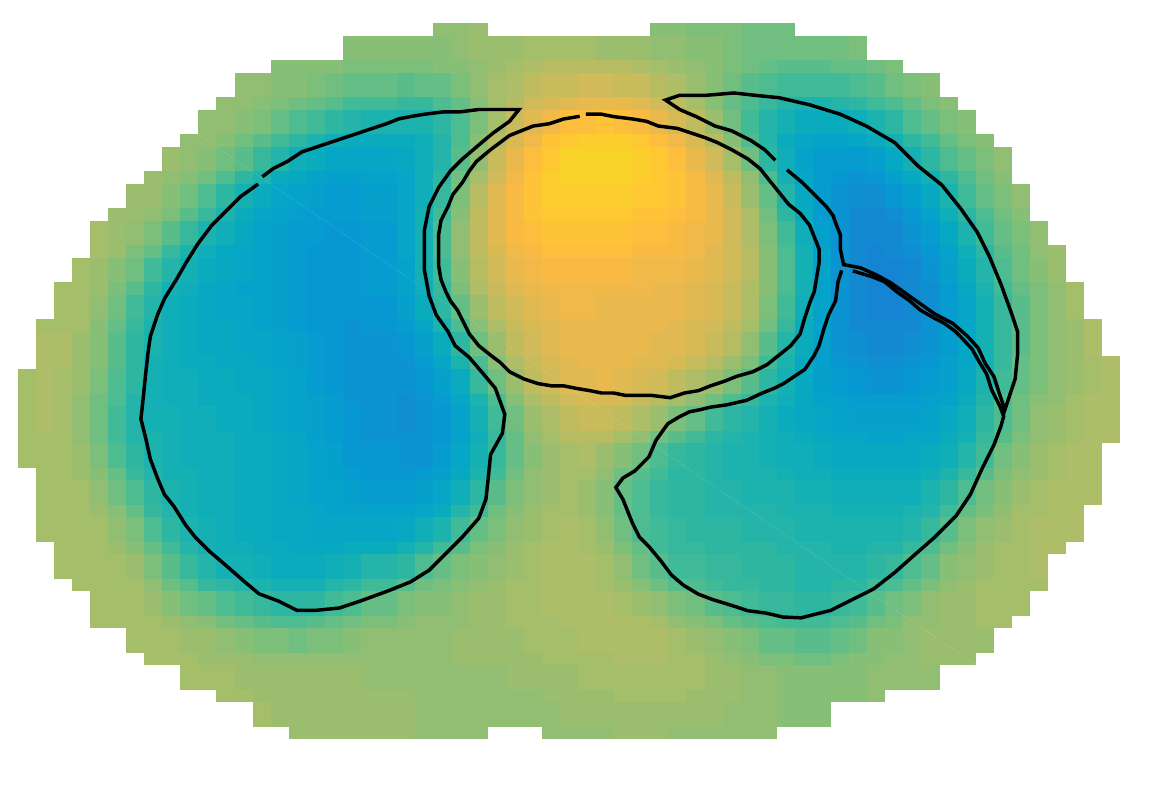}}\\

\vspace{1em}

\rotatebox{90}{\textbf{Permittivity}}
{\includegraphics[width=100pt]{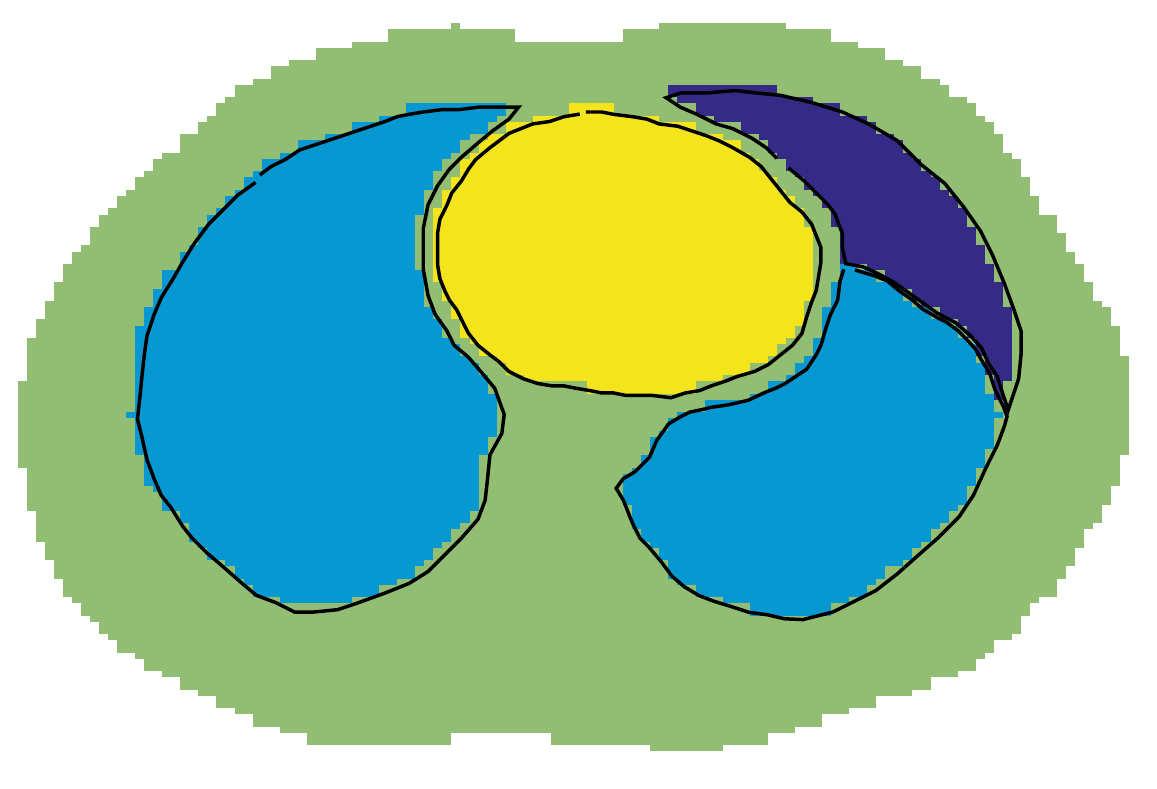}}
{\includegraphics[width=100pt]{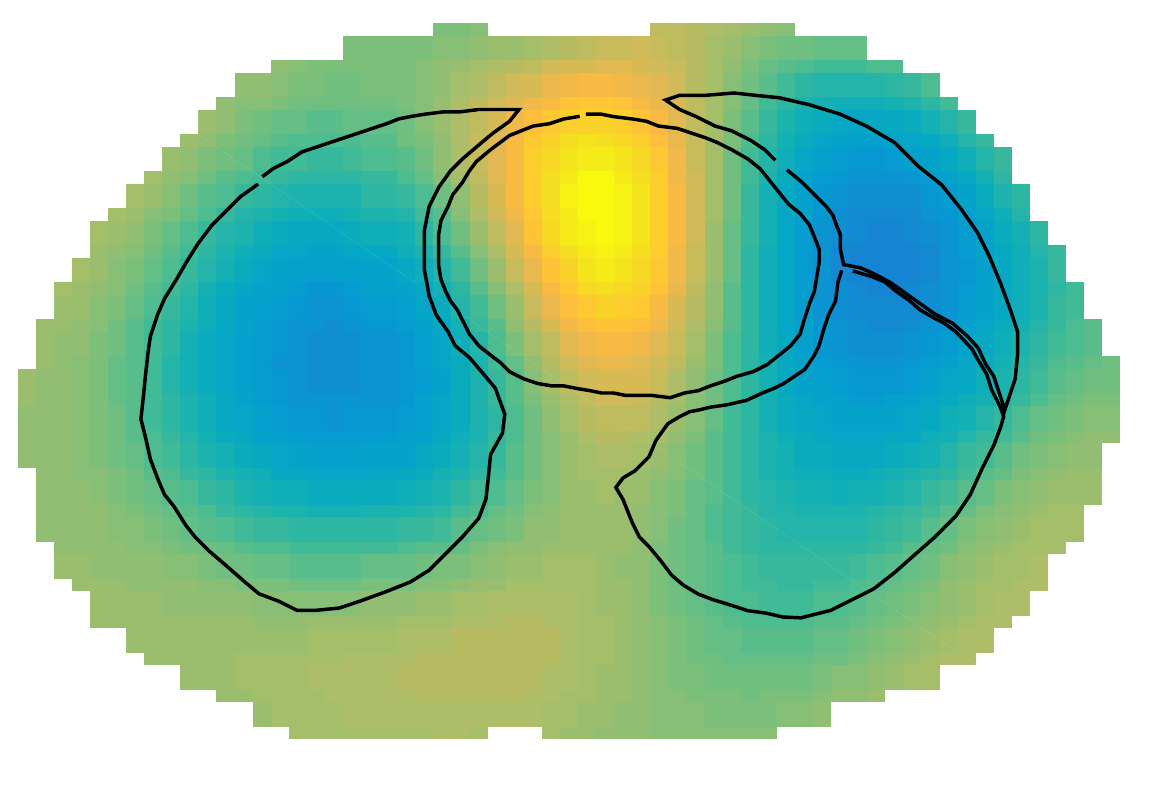}}
{\includegraphics[width=100pt]{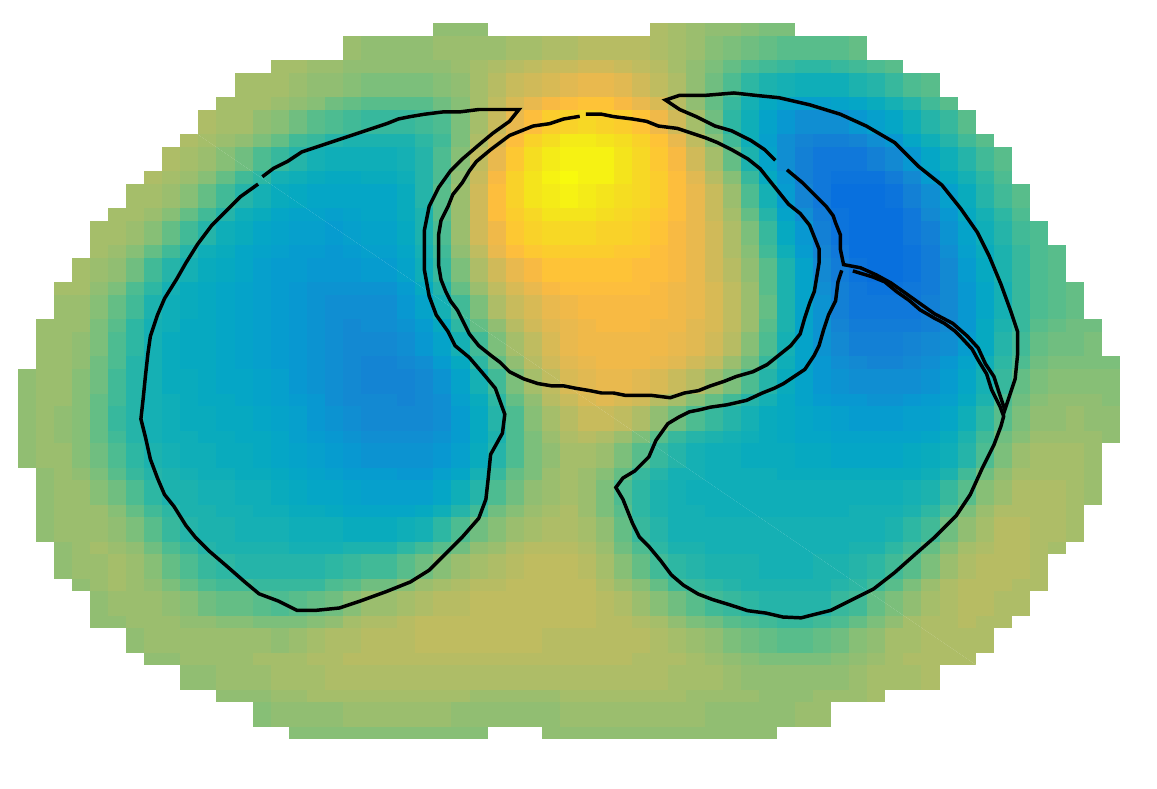}}\\

(a)\hspace{8em} (b) \hspace{8em} (c)
\caption{\label{fig:pneumo0o1_compare_strongest} Reconstructions for the pneumothorax example with 0.1\% noise plotted on the same scale.  Figure (a) is the true admittivity, (b) the initial D-bar reconstruction $\gDB$, and (c) the new admittivity $\gnew$ with $R_2=8.5$ and $\alpha=0$.}
\end{figure}

\begin{figure}[!h] 
\centering
{\textbf{Conductivity}}\hspace{7em}{\textbf{Permittivity}}\\
{\includegraphics[width=100pt]{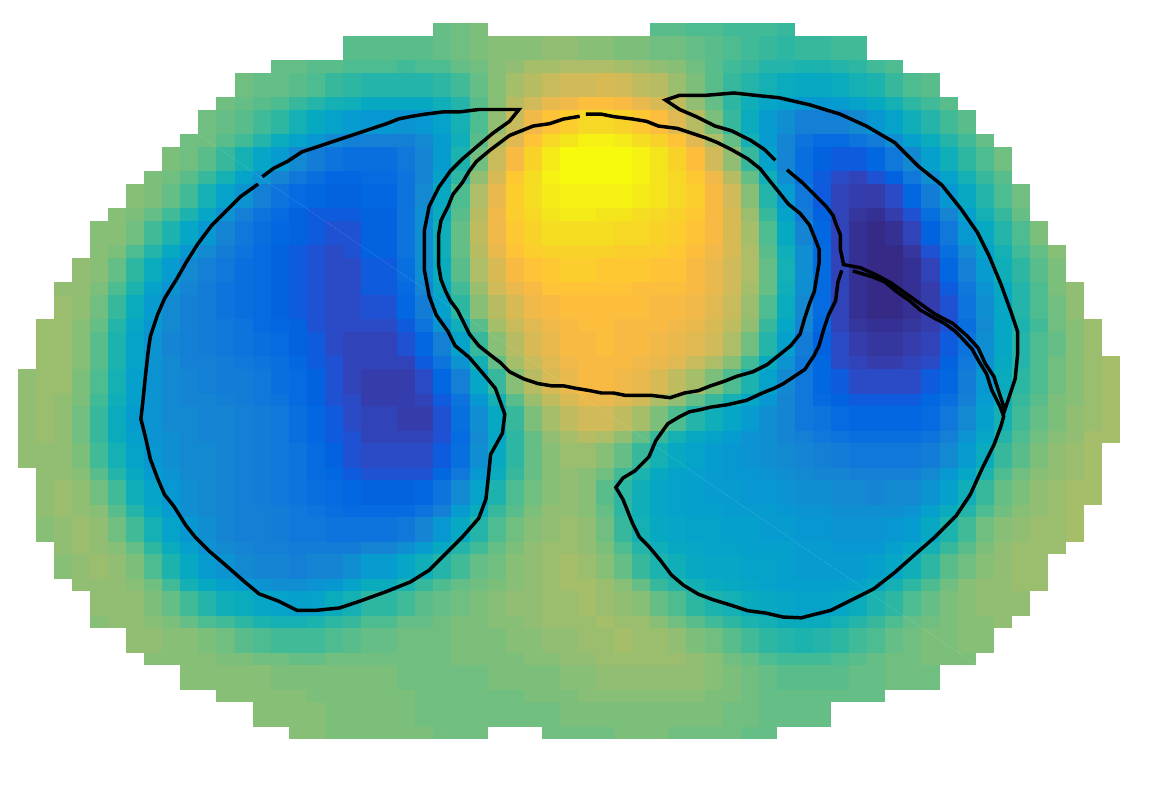}}\hspace{3em}
{\includegraphics[width=100pt]{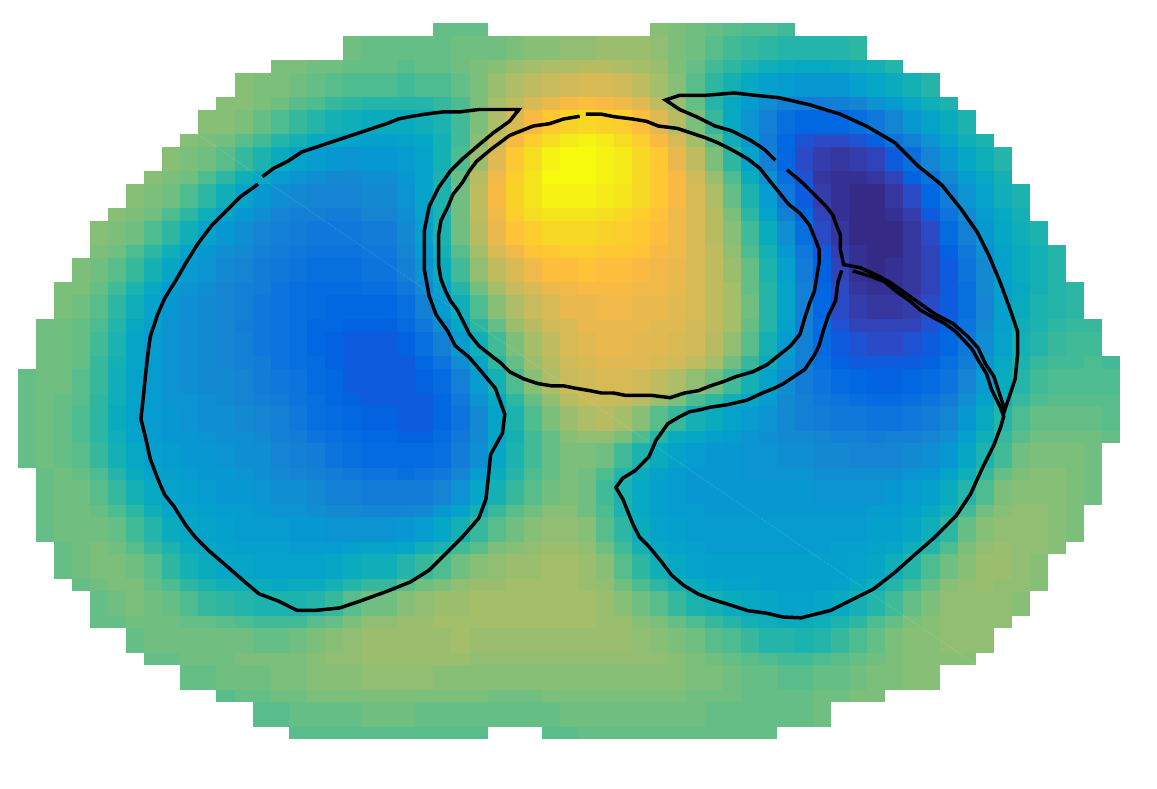}}
\caption{\label{fig:pneumo0o1_solo_strongest} Reconstruction $\gnew$ for the pneumothorax with 0.1\% noise with $R_2=8.5$ and $\alpha=0$.}
\end{figure}

\begin{figure}[!h] 

\centering
\rotatebox{90}{\textbf{Conductivity}}
{\includegraphics[width=100pt]{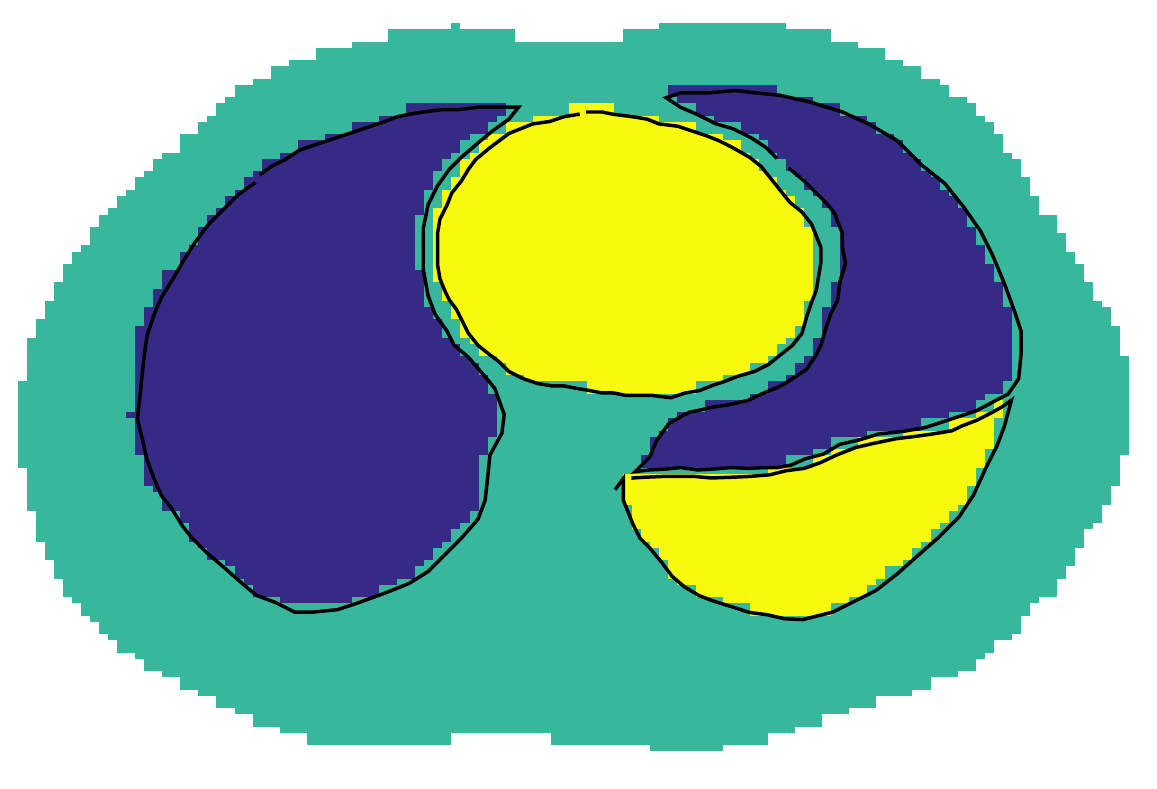}}
{\includegraphics[width=100pt]{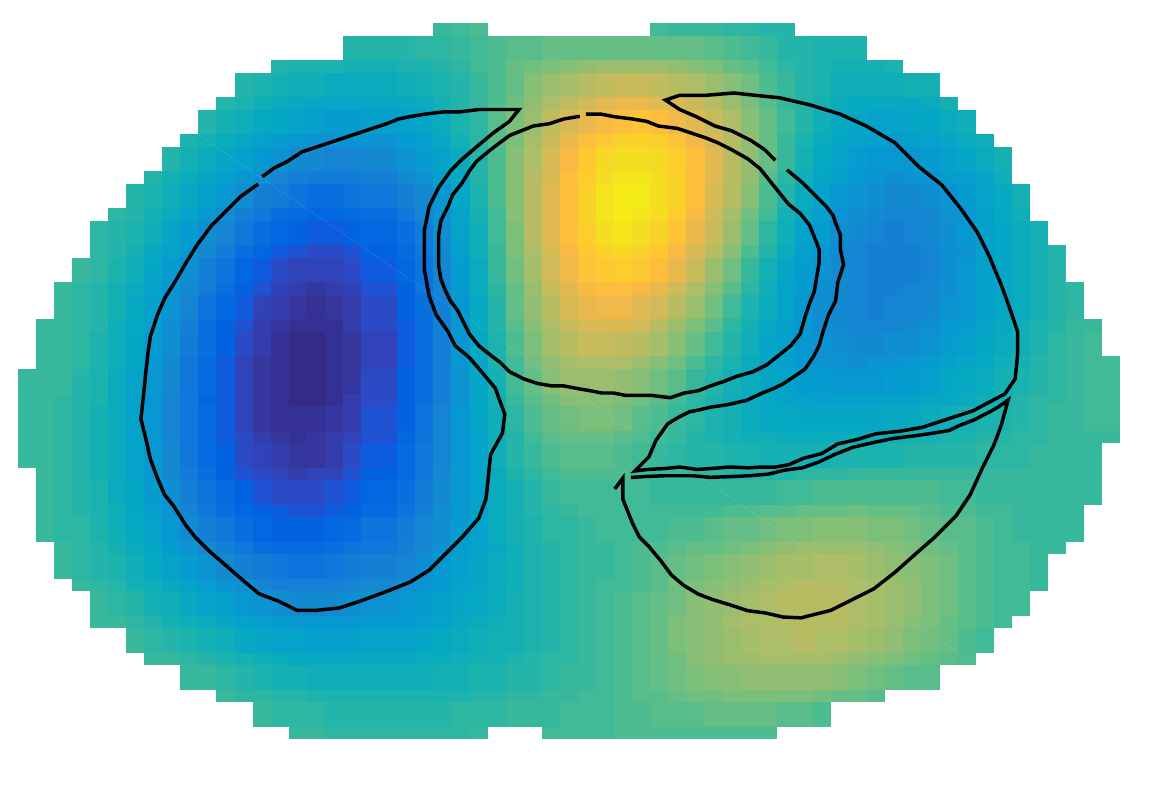}}
{\includegraphics[width=100pt]{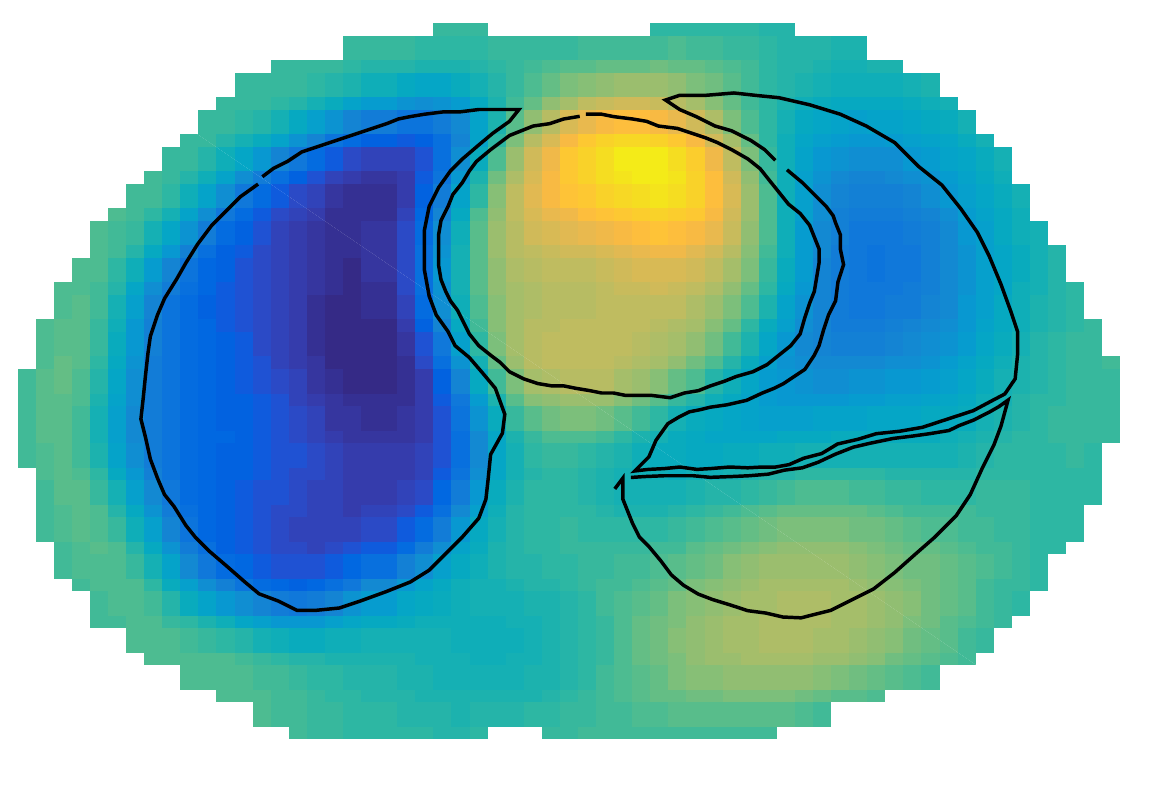}}\\

\vspace{1em}

\rotatebox{90}{\textbf{Permittivity}}
{\includegraphics[width=100pt]{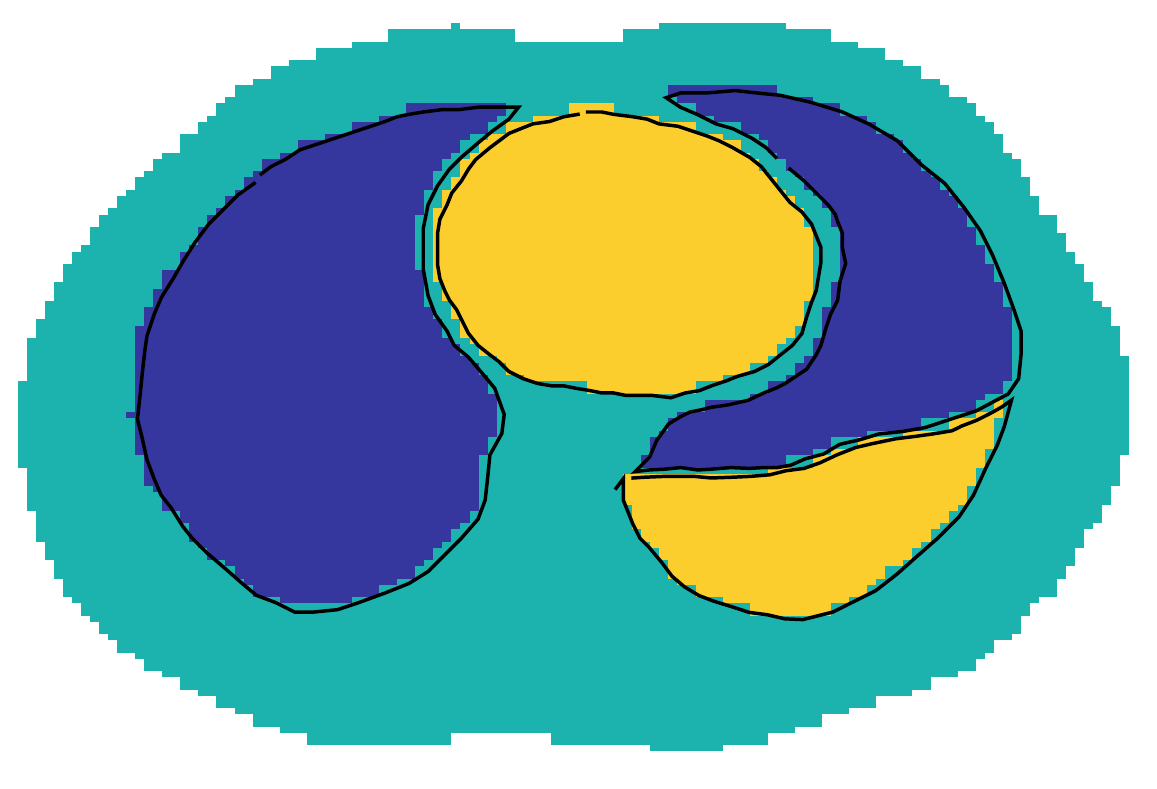}}
{\includegraphics[width=100pt]{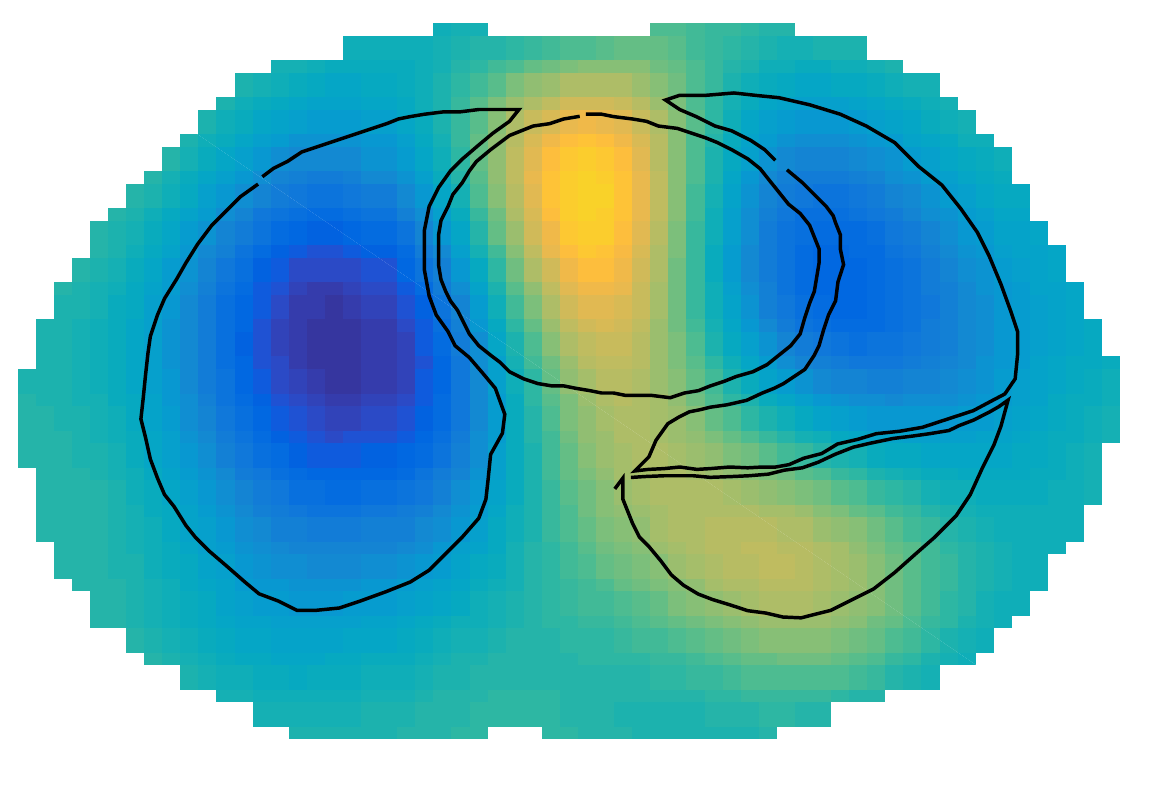}}
{\includegraphics[width=100pt]{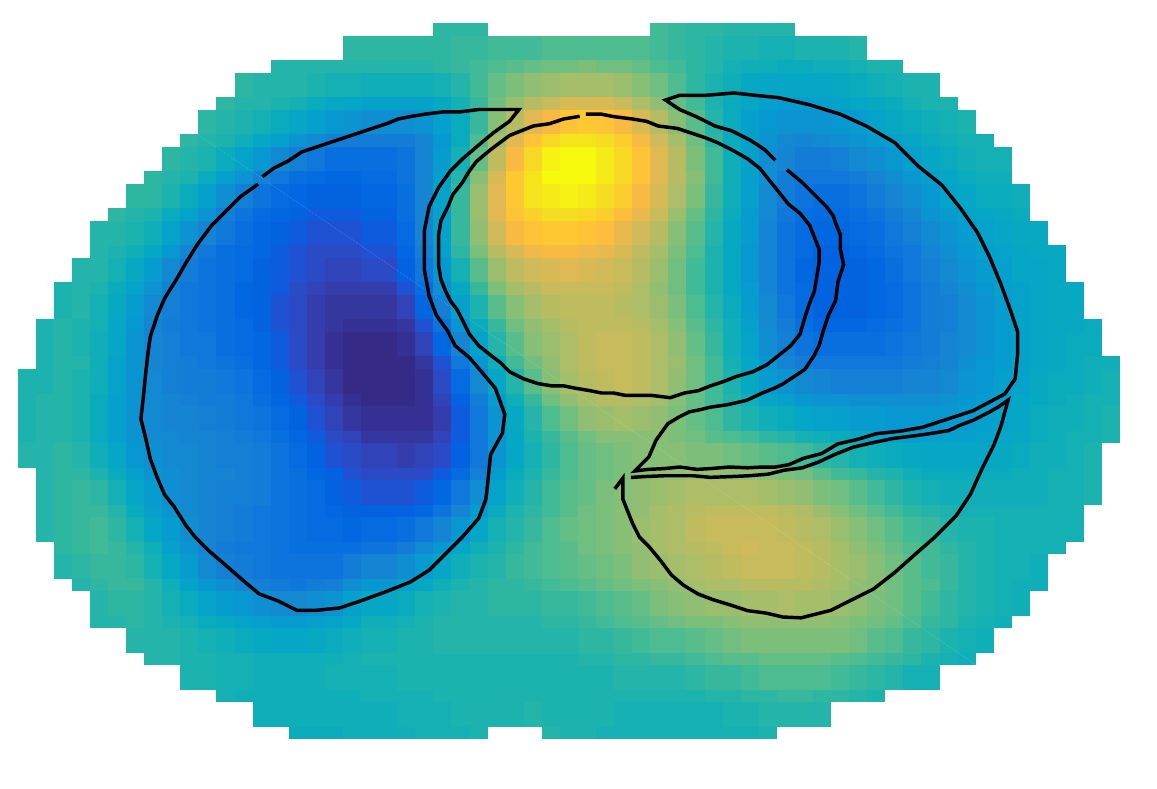}}\\

(a)\hspace{8em} (b) \hspace{8em} (c)
\caption{\label{fig:fluid0o1_compare_strongest} Reconstructions for the pleural effusion example with 0.1\% noise plotted on the same scale.  Figure (a) is the true admittivity, (b) the initial D-bar reconstruction $\gDB$, and (c) the new admittivity $\gnew$ with $R_2=11.0$ and $\alpha=0$.}
\end{figure}

\begin{figure}[!h] 
\centering
{\textbf{Conductivity}}\hspace{7em}{\textbf{Permittivity}}\\
{\includegraphics[width=100pt]{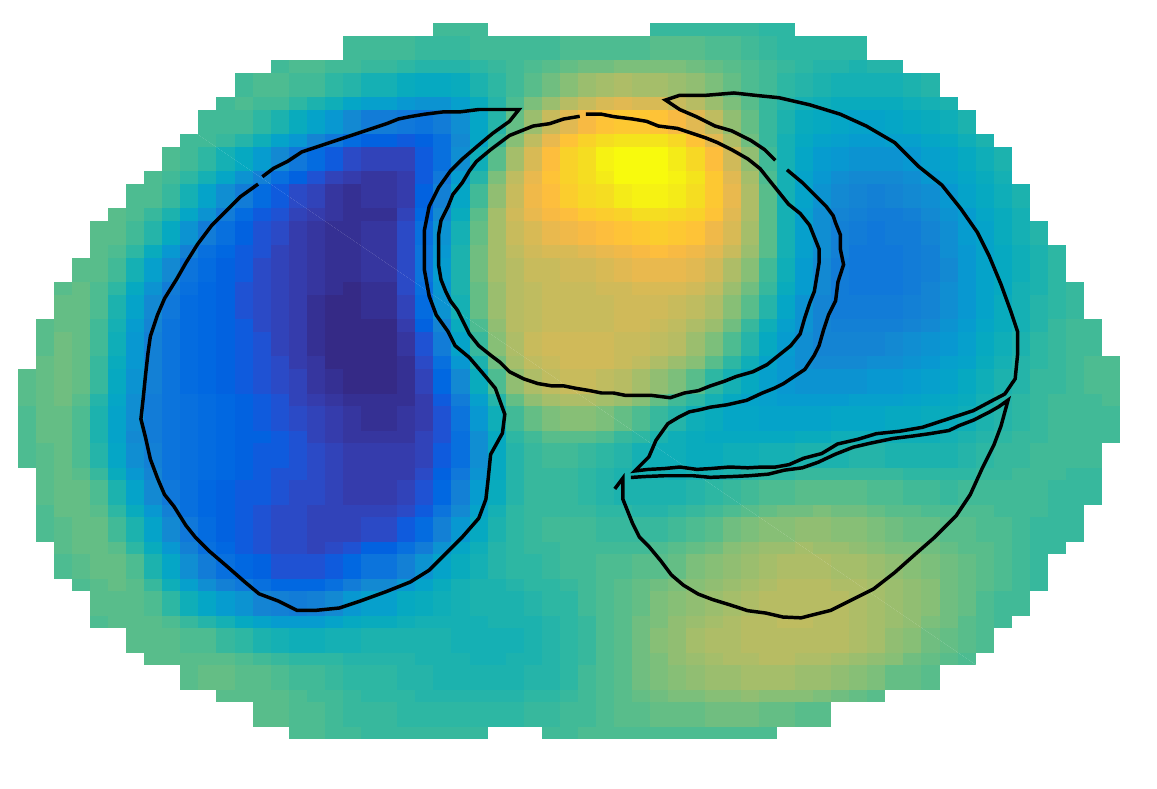}}\hspace{3em}
{\includegraphics[width=100pt]{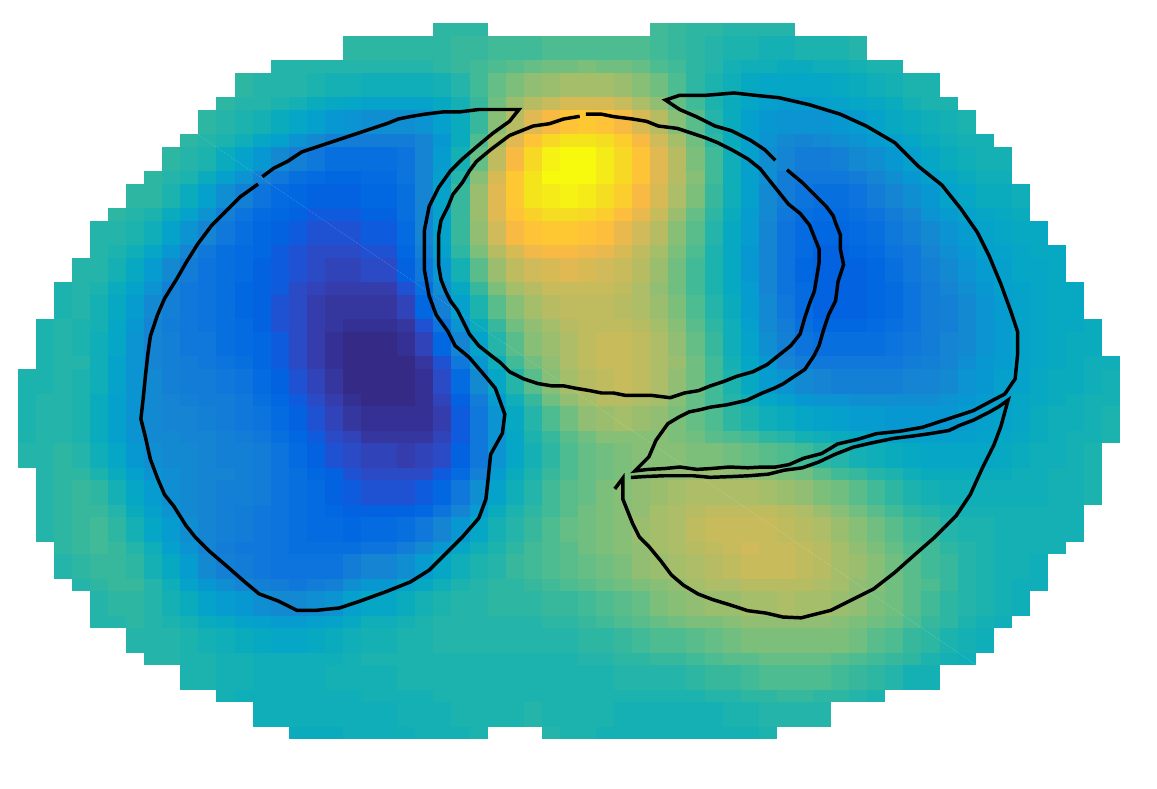}}
\caption{\label{fig:fluid0o1_solo_strongest} Reconstruction $\gnew$ for the pleural effusion with 0.1\% noise with $R_2=11$ and $\alpha=0$.}
\end{figure}

Figures~\ref{fig:pneumo0o1_compare_strongest} and \ref{fig:fluid0o1_compare_strongest} include side-by-side images of the (a) truth, (b) standard D-bar reconstruction with no prior, and (c) the  reconstructed conductivity and permittivity images with the strongest weights on the prior considered here, all displayed on the same scale for ease of comparison.  The true boundaries of the organs and pathologies are superimposed with black outlines.  Figures~\ref{fig:pneumo0o1_solo_strongest} and \ref{fig:fluid0o1_solo_strongest} show the new reconstructions alone for $\alpha=0$ and $R_2=8.5$ for the 0.1\% pneumothorax example, and $R_2=11$ for the 0.1\% pleural effusion example, to demonstrate the spatial improvement in the reconstructions.

It is clear from all of these images that this method is highly effective when organ boundaries are known with some confidence for improving the reconstructions without any bias of prior knowledge of the pathology.  The influence of various qualities of prior knowledge of the boundary and organ boundaries is left for future work, as are results from experimental data.  In practice, this high quality knowledge of organ boundaries corresponds to electrodes placed in the same plane as a CT scan slice.  This can be accomplished with careful use of fiducial markers, and averaging of several slices to account for the fact that EIT electrodes are typically much higher than a CT scan slice, resulting in an image that corresponds to a much thicker slice.  

\afterpage

\section*{Acknowledgments}
S.~J.~Hamilton was supported by the 2015 Summer Faculty Fellowship from Marquette University.  The project described was additionally supported by Award Number 1R21EB016869-01A1 from the National Institute Of Biomedical Imaging And Bioengineering.  The content is solely the responsibility of the authors and does not necessarily represent the official view of the National Institute Of Biomedical Imaging And Bioengineering or the National Institutes of Health.

\bibliographystyle{IEEEtran}
\bibliography{bibliographyRefs}

\end{document}